\documentclass[sn-mathphys-num]{sn-jnl}
\usepackage{graphicx}
\usepackage[numbers,square]{natbib}
\usepackage{amsmath,amsfonts,amssymb,mathrsfs,cancel,mathdots}
\usepackage{tcolorbox}

\usepackage{pgfplots, tikz}
\usetikzlibrary{babel,arrows.meta, patterns}
\pgfplotsset{compat=1.16}  
\usepgfplotslibrary{fillbetween}  
\usepackage{mathtools}
\usepackage{subcaption,graphicx,epstopdf}
\usepackage{amsthm}
\usepackage{amssymb}

\newcommand{\ren}{\mathbb{R}^n}
\newcommand{\re}{\mathbb{R}}
\newcommand{\propf}{$f: \ren \rightarrow \re \cup \{+\infty\}~$}

\newcommand{\epi}{{\rm epi}}
\newcommand{\dom}{{\rm dom}}
\newcommand{\blockcomment}[1]{}

\newtheorem{theorem}{Theorem}%  
\newtheorem{proposition}[theorem]{Proposition}% 

\newtheorem{example}{Example}%
\newtheorem{lemma}{Lemma}%
\newtheorem{corollary}{Corollary}
\newtheorem{remark}{Remark}

\newtheorem{definition}{Definition}%

\def\an#1{{\color{blue}#1}}

\begin{document}

\title[Article Title]{On existence of solutions to nonconvex minimization problems}

\author{\fnm{Rohan} \sur{Rele}}\email{rrele@asu.edu}

\author{\fnm{Angelia} \sur{Nedi\'c}}\email{anedich@asu.edu}

\affil[1]{\orgdiv{School of Electrical, Computer and Energy Engineering}, \orgname{Arizona State University}}

\abstract{
We provide a unified framework for a systematic analysis of 
the existence of solutions to general nonconvex problems, relying on asymptotic and retractive cones for functions and sets. Using this framework
we develop new necessary and sufficient conditions for the 
existence of solutions to a general problem of minimizing a proper closed function over a closed, possibly unbounded, set. Towards the result,
we introduce cones of retractive directions for a set and a function, establishing some basic properties for them. We also investigate the relationships between the cone of retractive directions of a function and the cone of level sets of the function.
Using the cones of retractive directions we provide necessary and sufficient conditions for the existence of solutions that require an asymptotically bounded decay of a function,  
and a relation between the cones of retractive directions of the constraint set and the asymptotic cone of the objective function. Finally we refine these conditions for more structured problems.
}

\keywords{Existence of Solutions, Nonconvex Problems, Optimization, Asymptotic Cone, Asymptotic Function}

%%\pacs[JEL Classification]{D8, H51}

%%\pacs[MSC Classification]{35A01, 65L10, 65L12, 65L20, 65L70}

\numberwithin{theorem}{section}
\numberwithin{definition}{section}
\numberwithin{lemma}{section}
\numberwithin{corollary}{section}
\numberwithin{example}{section}
\counterwithin{figure}{section}

\def\an#1{{\color{blue}#1}}

\maketitle

\section{Introduction}
\label{sec-intro}
Our interest is in investigating conditions for the existence of solutions to general nonconvex minimization problems by developing a simple unified framework relying on asymptotic cones of functions and sets.
The existence of solutions has been extensively studied starting with
seminal work~\cite{Frank_Wolfe_Theorem_56} showing that a quadratic function (bounded from below) attains its minimum on a polyhedral set.
The result has been extended to the problem with a quasi-convex objective function in~\cite{Luo_Zhang_1999}. The work in~\cite{Belousov_Klatte_2002} has established that a convex polynomial attains its solution on a region described by finitely many convex polynomial inequalities, which in turn generalized the result established in~\cite{Terlaky-LP} for convex quadratic functions. 

The inherent difficulty in establishing the existence of solutions is due to directions in unbounded constraint set along which the function may decrease.
A framework to address this developed in~\cite{Auslender_1997, baiocchi1998,auslender_teboulle_2002, penot2006}, relying on concepts such as the asymptotic cone and asymptotic function to show existence of solutions. Subsequently, in~\cite{bertsekas_tseng_2006}, this condition on asymptotic directions has been extended to introduce retractive directions and prove existence of solutions via a nonempty level-set intersection approach. The work in~\cite{ozdaglar_tseng_2006} has developed the existence results for problems where the constraint sets are given by functional inequalities. More recently, a solution existence result has been provided in~\cite{Hieu_2021} for a general polynomial objective and a closed constrained set under a certain regularity condition. 
Sufficient conditions for existence of solutions to convex problems and convex-concave saddle point problems over convex sets can be found in~\cite{bertsekas_nedic_ozdaglar_2013}.
Recently, work in~\cite{fakhar2023noncoercive} has used 
generalized asymptotic functions to obtain sufficient conditions for the existence of solutions to equilibrium problems in an infinite-dimensional space and, as a special result, provides a sufficient condition for 
the existence of solutions to a minimization problem.

 Similar to the works~\cite{Auslender_1997, baiocchi1998,auslender_teboulle_2002,ozdaglar_tseng_2006, penot2006} 
we make use of the asymptotic cones and asymptotic functions.
Unlike the aforementioned work, we introduce a concept of asymptotically bounded decay of a function over a given constraint set (with respect to a coercive function) and use it to establish necessary and sufficient conditions for the existence of solutions.
We introduce the cone of retractive directions for a set
 and a function and investigate their properties. While retractive directions for sets have been silently used to define asymptotically linear sets in~\cite{auslender_teboulle_2002}, and also used in~\cite{bertsekas_tseng_2006, penot2006} with a slightly different, but more general definition developed for a family of nested sets, the cone of such directions has not been studied.  
 Using these cones and asymptotically bounded decay property we establish new necessary and sufficient conditions for existence of solutions.
the properties of the cones of retractive directions for a set and a function.
More specifically, the novel aspects of the work include:

\noindent
(1)~Defining a class of functions that have asymptotically bounded decay with respect to a coercive function, which can be viewed as generalization of the coercivity and super-coercivity~\cite{Cambini-Carosi-Coercive,Lara_2022}. This class of functions includes proper closed convex functions and functions with Lipschitz continuous $p$-derivatives, with $p\ge0$, such as general multivariate polynomial functions (Section~\ref{sec-asb-decay}).

%(2)~The asymptotic cone of a set~\cite{auslender_teboulle_2002} and the asymptotic cone of a function, which extends such a notion for a proper convex function, as introduced in Definition~2.5.2  of~\cite{auslender_teboulle_2002}, to any proper function.

\noindent
(2)~Introducing a cone of retractive directions of a set and a function, and exploring their properties. Interestingly, while nonempty, such cones need not be closed, as shown by an example (Section~\ref{sec-retractive-dirset}). We also 
explore the relationships of a cone of retractive directions of a function and that of a level set of the function. 
Surprisingly, while asymptotic cone of a level set is always contained in the asymptotic cone of the function, this is not true for the cone of retractive directions. Even more surprisingly, for two distinct lower-level sets of a function, there is no particular inclusion relation for the cones of their retractive directions, even when the function is convex. 
Finally, we provide a closed form expression for the asymptotic function of a nonconvex polynomial and the cone of retractive directions of a convex polynomial (Section~\ref{sec-retractive-dirfun}).

\noindent

(3)~Establishing necessary and sufficient conditions for the existence of solutions based on asymptotically bounded decay property and a relation for cones of retractive directions of the objective function $f$ and the constraint set $X$ (Theorems~\ref{thm-compact-case}--\ref{main-result}, Theorem~\ref{thm-necessity}, Section~\ref{sec-main}).
The paper ~\cite{penot2006} is critical in showing the necessity of the conditions.
We refine our results for general nonconvex problems to a more structured problems, such as convex, nonconvex with functional constraints and polynomial constraints (Sections~\ref{sec-conv}--\ref{sec-nonconv}),  thus obtaining new results that are more general than those in~\cite{Belousov_Klatte_2002,ozdaglar_tseng_2006,bertsekas_tseng_2006}.

This paper is organized as follows: Section~\ref{sec-notation} provides the necessary background required for the subsequent development. Section~\ref{sec-problem} introduces the notion of asymptotically bounded decay of a function, the cones of retractive directions of a set and a function, and investigates their properties. 
Section~\ref{sec-main} presents our main results and their proofs for the problems where the constraint set is generic and not specified via inequalities.
Section~\ref{sec-conv} and Section~\ref{sec-nonconv} revisits our main results for the cases when the problem has more structure.

\section{Notation and Terminology}
\label{sec-notation}
We consider the space $\ren$
 equipped with the standard Euclidean norm $\|\cdotp\|$ unless otherwise stated. 

We deal with the extended real-valued functions that take values in the set $\re\cup\{+\infty\}$. For a function $f$, we use $\dom f$ to denote its effective domain, i.e., $\dom f = \{x \in \ren \mid f(x) < +\infty \}$. The epigraph of a function $f$ is denoted by $\epi f$, i.e., $\epi f = \{(x,c) \in \ren \times \re \mid f(x) \leq c \}$. For any $\gamma \in \re$,
    the lower-level set of a function $f$ is denoted by $L_\gamma(f)$, i.e., $L_\gamma(f) = \{x \in \ren \mid f(x) \leq \gamma\}.$
    
Our focus is on proper and closed functions, which are the functions with the epigraph $\epi f$ that is nonempty and closed set.

The function \propf is said to be coercive on a set $X$ if $$\lim_{ \substack{\|x\| \to +\infty \\ x\in X }} f(x) = +\infty. $$
The notion of a $p$-supercoercive, with $p \geq 1$, function has been introduced  in~\cite{Cambini-Carosi-Coercive, Lara_2022}, which states that:
    a proper function $f$ is $p$-supercoercive if
        $\liminf_{\|x\| \rightarrow \infty} \frac{f(x)}{\|x\|^p} > 0.$

\subsection{Asymptotic Cones}
\label{sec-asym-cones}

We turn our attention to unbounded sets whose behavior ``at infinity" is captured by their asymptotic cones (see~\cite{auslender_teboulle_2002}), which are also referred to as horizon cones \cite{RockWets98}.
To define these cones, we start with
the definition of a sequence converging in a direction. 

\begin{definition}\label{def-dir}
    A sequence $\{x_k\} \subset \ren$ is said to converge in the direction $d\in \ren$ if there exists a scalar sequence $\{t_k\} \subset \re$ with $t_k \rightarrow +\infty$ such that $\lim_{k\rightarrow\infty} t_k^{-1}x_k = d.$
\end{definition}

\begin{definition}[Definition 3.3 of~\cite{RockWets98}]
\label{asymp-cone}
    Let $X \subseteq \ren$ be a nonempty set. The asymptotic cone of $X$, denoted by $X_\infty$, is the set of vectors $d \in \ren$ that are limits in the directions of any sequence $\{x_k\} \subset X$ i.e., 
    $$ X_\infty = \left\{ d \in \ren \mid \exists\{x_k\} \subset X,\,\exists \{t_k\} \subset \re, t_k \to+\infty \ \text{\,such that }  \lim_{k\rightarrow\infty} \frac{x_k}{t_k} = d\right\}.$$
    For an empty set $X$, we have 
$X_\infty=\{0\}$. 
\end{definition}
The set $X_\infty$ is always nonempty closed cone (Proposition~2.1.1 in~\cite{auslender_teboulle_2002}). Asymptotic cones are illustrated in Figure~\ref{fig:asymp-cone-eq-recc}
and Figure~\ref{fig:asymp-cone-neq-recc}.
 We conclude this section with a result for the asymptotic cone of the intersection of sets 
 (Proposition~2.1.9 in~\cite{auslender_teboulle_2002}, or 3.9 Proposition in~\cite{RockWets98}).
\begin{proposition}\label{prop-set-intersect}
    Let $C_i \subseteq \ren$, $i \in \mathcal{I}$, where  $\mathcal{I}$ an arbitrary index set.
Then,
 \[(\cap_{i \in \mathcal{I}} C_i)_\infty \subseteq \cap_{i \in \mathcal{I}} (C_i)_\infty.\] 
    The inclusion holds as an equality for closed convex sets $C_i$. 
\end{proposition}

\begin{figure}[hb!]
\begin{minipage}{.45\textwidth}
\center

      \begin{tikzpicture}[scale = 0.6]
  \begin{axis}[
    xlabel={$x$},
    ylabel={$y$},
    xmin=-2, xmax=2,
    ymin=-4, ymax=4,
    domain=-2:2,
    samples=100,
    axis lines=middle,
    axis line style={-}, 
    ticks=none, 
    clip=true, 
    ]

    \addplot [name path = positive, black, domain=-2:2] {x^2};
    \addplot [name path = top, teal!30, forget plot] coordinates {(-2,4) (2,4)};

    \addplot [name path = negative, black, domain=-2:2] {-x^2};
    \addplot [name path = bottom, teal!30, forget plot] coordinates {(-2,-4) (2,-4)};

    \addplot [teal!30] fill between [of = positive and top];
    \addplot [teal!30] fill between [of = negative and bottom];

    \node[circle,fill,red, inner sep=2pt] at (axis cs:0,0) {};

    \addplot [name path = cone, red, domain = -4:4, line width = 1.0pt, <->] coordinates {(0,4) (0,-4)};
  \end{axis}
  \end{tikzpicture}
  \captionof{figure}{Set $X=\{(x,y)\in\re^2\mid x^2\le |y| \}$ and its asymptotic cone $X_\infty=\{(0,y)\mid y\in\re\}$.}
  \label{fig:asymp-cone-eq-recc}
 \end{minipage}%
 \qquad
 \begin{minipage}{.45\textwidth}
 \center
  \begin{tikzpicture}[scale = 0.6]
     \begin{axis}[
    xlabel={$x$},
    ylabel={$y$},
    xmin=-2, xmax=2,
    ymin=-4, ymax=4,
    domain=-2:2,
    samples=100,
    axis lines=middle,
    axis line style={-}, % Removes arrow heads from axis lines
    ticks=none, % Removes ticks from the axes
    clip=true, % Prevents clipping of the plot
    ]

    \addplot [name path = positive, black, domain=-2:2] {x^2};
    \addplot [name path = bottom, teal!30, forget plot] coordinates {(-2,-4) (2,-4)};

    \addplot [teal!30] fill between [of = positive and bottom];
    \addplot [name path = rcone, black, domain = -4:0, line width = 0.5pt] coordinates {(0,0) (0,-4)};
    \addplot [name path = x-axis, black, domain = -2:2, line width = 0.5pt] coordinates {(-2,0) (2,0)};
    \end{axis}
    
  \end{tikzpicture}
  \captionof{figure}{Set $X=\{(x,y)\in\re^2\mid x^2\ge y \}$ and its asymptotic cone $X_\infty=\re^2$.}
  \label{fig:asymp-cone-neq-recc}
 \end{minipage}
\end{figure}

\subsection{Asymptotic Functions}
Applying the concept of the asymptotic cone to the epigraph of a proper function leads the notion of an asymptotic function, which is also referred to as a horizon function\cite{RockWets98}. We elect to define it according to \cite{auslender_teboulle_2002}, as stated below. 

\begin{definition}\label{asymp-def}
    For any proper function \propf, there exists a unique function $f_\infty: \ren \rightarrow \re \cup \{-\infty, +\infty\}$ associated with $f$ such that $\epi f_\infty = (\epi f)_\infty$. The function $f_\infty$ is said to be the asymptotic function of $f$.
\end{definition}

A useful analytic representation of an asymptotic function $f_\infty$ was originally obtained in~\cite{Auslender_1997} and, also, given in Theorem~2.5.1 of~\cite{auslender_teboulle_2002}.

\begin{theorem}\label{thm-asfun-rep}
    For any proper function \propf the asymptotic function $f_\infty$ is given by 
    \begin{equation}\label{AF1}
         f_\infty(d) = \liminf_{\substack{d' \rightarrow d \\ t \rightarrow +\infty}} \frac{f(td')}{t},
    \end{equation}
    or equivalently,
    \begin{equation}\label{AF2}
        f_\infty(d) = \inf \left\{ \liminf_{k \rightarrow \infty} \frac{f(t_kd_k)}{t_k} ~ \bigg| ~ t_k \rightarrow +\infty, d_k \rightarrow d \right\}, 
    \end{equation}
    where the infimum is taken over all sequences $\{d_k\}\subset \ren$ and $\{t_k\} \subset \re$.
\end{theorem}

An asymptotic function has some basic properties inherent from its definition, namely that $f_\infty$ is closed and positively homogeneous function since its epigraph is the closed cone $(\epi f)_\infty$.

Further, the value $f_\infty(0)$ is either finite or $f_\infty(0) = -\infty$. If $f_\infty(0)$ is finite, then it must be that $f_\infty(0) = 0$ by the positive homogeneity property. As a consequence, for a proper function $f$ we have that 
$$0\in \{d\mid f_\infty(d)\le0\},$$ implying that 
$\{d\mid f_\infty(d)\le0\}\ne\emptyset.$
The directions $d$ such that $f_\infty(d)\le0$ will be particularly important in our subsequent development.

To this end, we will term them as asymptotic directions of a function, and use these directions to define the asymptotic cone of a function, as follows.

\begin{definition}\label{def-rec-f}
For a proper function $f$, we say that a direction $d$ is an asymptotic direction of $f$ if $f_\infty(d)\le 0$. 
The asymptotic cone of $f$, denoted by $\mathcal{K}(f)$, is the set of all asymptotic directions of $f$, 
i.e.,
\[\mathcal{K}(f) =\{d\mid f_\infty(d)\le 0\}.\]
\end{definition}
An asymptotic direction of a function has been given in Definition~3.1.2 of~\cite{auslender_teboulle_2002},  while the asymptotic cone of a function has been defined for a proper convex function in Definition~2.5.2 of~\cite{auslender_teboulle_2002}. Here, we adopt the same definition for any proper function.

We next provide a key result for the asymptotic cones of lower-level sets of a proper function (Proposition~2.5.3 of
\cite{auslender_teboulle_2002} and 3.23~Proposition of~\cite{RockWets98}). 

\begin{proposition}

\label{prop-level-set}
    For a proper function $f$ and $\alpha \in \re$, we have $\left(L_\alpha(f)\right)_\infty \subseteq L_0(f_\infty)$ i.e.,
    \begin{equation*}%\label{level-set-inclusion}
        \{x \mid f(x) \leq \alpha\}_\infty \subseteq 
         \mathcal{K}(f).
    \end{equation*}
    The inclusion is an equality when $f$ is proper, closed, and convex.
\end{proposition}

Finally, we state an existence of solutions result that will also be used in the sequel.

\begin{proposition}
\label{prop-exist-basic}
Let $X\subseteq\mathbb{R}^n$ be a closed set, and let \propf be a proper closed function with $X\cap\dom f\ne\emptyset$. If $f$ is coercive over $X$, i.e.,
\[\liminf_{\|x\|\to\infty\atop x\in X} f(x)=+\infty,\]
then the problem $\inf_{x\in X} f(x)$ has a finite optimal value and an optimal solution exists.
\end{proposition}
\begin{proof}
The result follows by applying Theorem~2.14 of~\cite{beck-book} to the function $f+\delta_X$, where $\delta_X$ is the characteristic function of the set $X$.

\end{proof}

\section{Asymptotically Bounded Decay, Cones of Retractive Directions}\label{sec-problem}

Here, we introduce a notion of asymptotically bounded decay of a function (Subsection~\ref{sec-asb-decay}), and define the cones of retractive directions for sets and functions and study their properties (Subsections~\ref{sec-retractive-dirset} and~\ref{sec-retractive-dirfun}, respectively).

\subsection{Asymptotically Bounded Decay}
\label{sec-asb-decay}

We introduce a condition that prohibits the decrease of $f$ to infinity at a certain rate, which can be viewed as a generalization of the coercivity property.

\begin{definition}\label{def-asymp-decay}
    A proper function \propf is said to exhibit asymptotically bounded decay with respect to a proper function $g: \ren \to \re \cup \{+\infty\}$ on a set $X \subseteq \ren$ if 
    \begin{equation}\label{decay-condition}
        \liminf_{\substack{\|x\| \rightarrow \infty \\ x \in X}} \frac{f(x)}{g(x)} > -\infty.
    \end{equation}
    We say a function exhibits asymptotically bounded decay with respect to $g$ if $X = \ren$. 
\end{definition}

The condition in~\eqref{decay-condition} prohibits the function $f$ from
approaching $-\infty$, along the points in the set $X$, faster than the function $g$. 
Note that, if $f$ is $p$-supercoercive on $X$, with $p \geq 1$, then the bounded decay condition in~\eqref{decay-condition} is satisfied with $g(x)=\|x\|^p$.
The class of functions that have this asymptotically bounded decay is rather rich, as illustrated in the following examples.

\begin{example}\label{proper-finite}
    Let $X \subseteq \ren$ be a nonempty set. Let \propf be a proper function with a finite minimum on $X$, i.e.,
    $f^* = \inf_{x \in X} f(x)>-\infty$.
    Then, we have
    $f(x)\ge f^*$ for all $x\in  X,$
    implying that
    \begin{equation*}
        \liminf_{\substack{\|x\| \to \infty \\ x \in X}} \frac{f(x)}{\|x\|}\geq \liminf_{\substack{\|x\| \to \infty \\ x \in X}} \frac{f^*}{\|x\|}=0 .
    \end{equation*}
    Thus, if $\inf_{x \in X} f(x)>-\infty$, then $f$ satisfies Definition \ref{def-asymp-decay}
    with $g(x)=\|x\|$.
    \qed
\end{example}

Next, we show that any proper convex function
satisfies Definition \ref{def-asymp-decay}
    with $g(x)=\|x\|$ due to the linear underestimation property of a convex function.

\begin{example}\label{ex-convex}
Let $f$ be a proper convex function such that $X\cap\dom f\ne \emptyset$, and let $x_0$ be a point in the relative interior of $\dom f$. Then, by Theorem 23.4 of~\cite{rockafellar-1970a}, the subdifferential set $\partial f(x_0)$ is nonempty. Thus, by the convexity of $f$ we have for a subgradient $s_0$ of $f$ at the point $x_0$ and for all $x\in X$,
\begin{align*}
  f(x) \ge f(x_0)+ \langle s_0,x-x_0\rangle  \ge f(x_0)-\|s_0\| \|x-x_0\|,
\end{align*}
implying that 
\[ \liminf_{\|x\|\to\infty\atop x\in X}\frac{f(x)}{\|x\|}\ge \liminf_{\|x\|\to\infty\atop x\in X}\frac{f(x_0)-\|s_0\| \|x-x_0\|}{\|x\|}=-\|s_0\|. \]
\qed
\end{example}

Our next example shows that if a function $f$ satisfies 
$f_\infty(d)\ge0$ for all nonzero $d\in\re^n$, then
$f$ satisfies Definition \ref{def-asymp-decay}
    with $g(x)=\|x\|$.
    %exhibits asymptotically bounded decay with respect to $\|x\|$.
\begin{example}\label{ex-finfty} Let $X$ be a nonempty set and $f$ be a proper function such that $X\cap\dom f\ne \emptyset$. If
$f_\infty(d)\ge0$ for all nonzero $d\in X_\infty$, then
 \[
            \liminf_{\|x\| \to \infty\atop x\in X} \frac{f(x)}{\|x\|} = 
            \liminf_{\|x\| \to \infty} \frac{f(x)+\delta_X(x)}{\|x\|} =
            \liminf_{\substack{d = x\|x\|^{-1} \\ \|x\| \to \infty}}\frac{f(\|x\|d)+\delta_X(\|x\|d)}{\|x\|} , 
            \]
            where
            $\delta_X(x)=0$ when $x\in X$ and otherwise, $\delta_X(x)=+\infty$.
            Therefore,
    \begin{align*}
    \liminf_{\|x\| \to \infty\atop x\in X} \frac{f(x)}{\|x\|} & \ge 
            \liminf_{d'\to d\atop t \to \infty} \frac{f(td') +\delta_X(td')}{t}  \cr
            &\ge f_\infty(d)+(\delta_X)_\infty(d).\end{align*}
Since ($\delta_X)_\infty(d) = \delta_{X_\infty}(d)$ and $f_\infty(d)\ge0$ for all nonzero $d\in X_\infty$, it follows that 
        $$ \liminf_{\|x\| \to \infty \atop x\in X} \frac{f(x)}{\|x\|} \ge 0, $$
       and $f$ exhibits asymptotically bounded decay on $X$ with respect to $\|x\|$.
       \qed
       \end{example}
Finally, we show that a $p$-times differentiable function, with Lipschitz continuous $p$th derivatives, satisfies Definition \ref{def-asymp-decay}
   with $g(x)=\|x\|^{p+1}$. 

\begin{example}\label{ex-lip}
    Consider a function \propf with Lipschitz continuous $p$th derivatives on an open convex set containing the set $X$, with $p\ge0$. When $p=0$, the function is simply Lipschitz continuous. Let the $p$th
    derivative have a Lipschitz constant $L_p>0$, i.e.,
        $$\|D^p f(x) - D^p f(x') \| \leq L_p\|x-x'\|
        \qquad \hbox{for all $x, x' \in \dom f$},$$
        where $D^pf(x)$ denotes the $p$th derivative of $f$ at a point $x$.
    Then, by Equation (1.5) of~\cite{nesterov2021}, we have that 
    \[|f(x)-\Phi_{x_0,p}(x)|\le \frac{L_p}{(p+1)!}\|x-x_0\|^{p+1}
    \qquad \hbox{for all $x, x_0 \in \dom f$},\]
    where $\Phi_{x_0,p}(x)$ is the $p$th order Taylor approximation of $f$ at the point $x_0$, i.e.,
    \[\Phi_{x_0,p}(x)=\sum_{i=1}^p\frac{1}{i!} D^if(x_0)[x-x_0]^i,\]
    with $[h]^i$ denoting the vector consisting of $i$ copies of a vector $h$, and with $[h]^0=1$ when $i=0$.
    Then, for $x_0\in \dom f$ arbitrary but fixed, we have that
    \[f(x)\ge \Phi_{x_0,p}(x) - \frac{L_p}{(p+1)!}\|x-x_0\|^{p+1}
    \qquad \hbox{for all $x\in \dom f$},\]
    implying that 
    \begin{align*}\liminf_{k\to\infty\atop x\in X} 
    \frac{f(x)}{\|x\|^{p+1}}
    \ge \liminf_{k\to\infty\atop x\in X}\left\{\frac{\Phi_{x_0,p}(x)}{\|x\|^{p+1}}- \frac{L_p}{(p+1)!}\frac{\|x-x_0\|^{p+1}}{\|x\|^{p+1}}\right\}=
    -\frac{L_p}{(p+1)!}, 
    \end{align*}
    where we used the fact that $\lim_{\|x\|\to\infty} \Phi_{x_0,p}(x)/\|x\|^{p+1}=0$.
   
    Hence, $f$ exhibits asymptotically bounded decay 
    on $X$ with respect to $g(x) = \|x\|^{p+1}$.
    \qed
\end{example}
As a special case of Example~\ref{ex-lip}, a function with Lipschitz continuous gradients 
exhibits asymptotically bounded decay with respect to the function $g(x)=\|x\|^{2}$.

Moreover, note that multivariate polynomials are a special class of functions that fall under Example~\ref{ex-lip}. In particular,
a multivariate polynomial of order $m$, with $m\ge1$, has a constant $m$th derivative, say equal to $B$. Thus, the $(m-1)$st derivative is Lipshitz continuous with the constant $B$.
By Example~\ref{ex-lip}, a multivariate polynomial
has asymptotically bounded decay with respect to the function $g(x)=\|x\|^{m}$.

\subsection{Cone of Retractive Directions of Sets}
\label{sec-retractive-dirset}

The key notion that we use throughout the rest of this paper is that of a retractive direction. For a nonempty set, a retractive direction is defined as follows.

\begin{definition}\label{retractive-def}
Given a set $X$, 
    a direction $d \in X_\infty$ is said to be retractive direction of $X$ if for any sequence $\{x_k\} \subseteq X$ converging in the direction $d$ and for any $\rho > 0$, there exists an index $K$ (depending on $\rho$) such that  
    \begin{equation}\label{retractive-eq}
        x_k - \rho d \in X \quad \text{ for all } k \geq K.
    \end{equation} The set of retractive directions of a set $X$ is denoted by $\mathcal{R}(X)$. We say that the set $X$ is retractive if $\mathcal{R}(X) = X_\infty$. 
\end{definition}

Note that, by definition, we have $\mathcal{R}(X) \subseteq X_\infty$, $0 \in \mathcal{R}(X)$ and, for $X$ empty, $\mathcal{R}(X) =X_\infty=\{0\}$.
We next provide examples of a convex set and a nonconvex set that have no nonzero retractive direction.  

\begin{example}\label{nonretractive-convex-ex}
    Consider the epigraph of the scalar function $f(s) = s^2$ i.e., $X = \{(s,\gamma) \in \re^2\mid s^2 \leq \gamma \}.$ Let $\{x_k\} \subseteq X$ be given by $x_k = (\sqrt{k}, k). $ Then, $\|x_k\| \to\infty$ as $k \rightarrow \infty$. For any $\lambda>0$, we have with $t_k=\|x_k\|/\lambda$,
    $$ \lim_{k\to\infty} \frac{x_k}{t_k}
    =\lambda
    \lim_{k\to\infty} \frac{x_k}{\|x_k\|} 
    =\lambda \lim_{k\to\infty}\left( \frac{1}{\sqrt{k + 1}}, \frac{\sqrt{k}}{\sqrt{k+1}} \right)=(0,\lambda).$$ 
    Thus, $(0,\lambda) \in X_\infty$ for any $\lambda>0$ and, furthermore, 
    $$ x_k - (0,\lambda)= (\sqrt{k}, k - \lambda)\not\in X\qquad\hbox{for all }k\ge 1.$$ 
    Thus, $d=(0,\lambda)$ is not a retractive direction of $X$ for any $\lambda>0$, and
    $\mathcal{R}(X)=\{0\}$. 
    \qed
\end{example}

\begin{example}\label{nonretractive-nonconvex-ex}
Consider the set $X=\{(x_1,x_2)\in\re^2\mid x_1^2\le |x_2|\}$ (see Fig.~\ref{fig:asymp-cone-eq-recc}). Similar to Example~\ref{nonretractive-convex-ex}, we can see that the directions $(0,\lambda)$ and $(0,-\lambda)$ are not retractive directions of $X$ for any $\lambda>0$. Hence, $\mathcal{R}(X)=\{0\}$.
\qed
\end{example}

Now we highlight some related definitions. Most notable is Definition~2.3.1 in~\cite{auslender_teboulle_2002} which 
defines an \textit{asymptotically linear set}, as follows: a closed set $X \subseteq \ren$ is said to be asymptotically linear if for every $\rho > 0$ and each sequence $\{x_k\} \subseteq X$ that satisfies $x_k \in X$, $\|x_k\| \rightarrow +\infty$ and $x_k\|x_k\|^{-1} \rightarrow \overline{x}$, there exists an index  $K$ such that $x_k - \rho\overline{x} \in X$ for all $k \geq K$. Note that the directions $\bar x$ involved in this definition have unit norm,
and such directions are retractive according to our Definition~\ref{retractive-def}. Further, since $\bar x$ cannot be zero, the set of such directions is a subset of $\mathcal{R}(X)$. 

By the definition of an asymptotically linear set, 
a polyhedral set $X$ is retractive i.e., $\mathcal{R}(X)=X_\infty$, as detailed below.

\begin{example}\label{ex-polyhedral-retract}
    Consider a polyhedral set $X$.

    By Proposition 2.3.1 in \cite{auslender_teboulle_2002}, an asymptotically polyhedral set $X \subseteq \ren$ is asymptotically linear. Since the simplest case of an asymptotically polyhedral set is a polyhedral set, it follows that $X$ is asymptotically linear, i.e.,
    $\mathcal{R}(X)=X_\infty$. 
    \qed
\end{example}

Another related definition is Definition 1 in~\cite{bertsekas_tseng_2006}, which considers
$x_k \in C_k$ for an infinite sequence of nested sets $\{C_k\} \subseteq \ren$ i.e., $C_{k+1} \subseteq C_k$ for all $k$. The directions $d$ of interest are obtained in the limit, as follows:
\[\lim_{k\to\infty} \frac{x_k}{\|x_k\|}=\frac{d}{\|d\|}.\]
By letting $C_k = X$ for all $k$, 
a direction $d$ is retractive according to Definition~1 in~\cite{bertsekas_tseng_2006}
if, for any associated sequence $\{x_k\} \subseteq X$ with $\|x_k\| \rightarrow$ and $x_k/\|x_k\|\to d/\|d\|$, we have that $x_k - d \in X$ for all sufficiently large $k$.
According to this definition, given a retractive direction $d$ and its associated sequence $\{x_k\}$, for any $\rho>0$, we have that 
\[\lim_{k\to\infty} \frac{x_k}{\|x_k\|}=\frac{d}{\|d\|}=
\frac{\rho d}{\|\rho d\|},\]
implying that $\{x_k\}$ is also associated sequence for the direction $\rho d$ for any $\rho>0$.
Thus, the condition $x_k-d\in X$ for all $k$ large enough can be written 
as $
x_k-\rho d\in X$ for all large enough $k$, implying that
$d$ is a
retractive direction according to Definition~1 in~\cite{bertsekas_tseng_2006} and, also, according to our Definition~\ref{retractive-def} with $t_k = \|x_k\|$ for all $k$. 

Consider now a sequence $\{x_k\}\subset X$ converging in a nonzero direction $d$, i.e., for some scalar sequence $\{t_k\}$ with $t_k\to\infty$, we have
\[\lim_{k\to\infty}\frac{x_k}{t_k}=d\qquad\hbox{ with }d\ne0.\]
Then,
\[\lim_{k\to\infty}\frac{x_k}{\|x_k\|}=\lim_{k\to\infty} \frac{x_k/t_k}{\|x_k\|/t_k}=\frac{d}{\|d\|}. \]
If $d$ is retractive according to our Definition \ref{retractive-def},
then by letting $\rho=\|d\|$, we conclude the $d$ is also retractive according to Definition~1 in~\cite{bertsekas_tseng_2006}. As stated previously, the only direction $d$ for which the equivalence does not hold is $d = 0$. That is, $0 \in \mathcal{R}(X)$ by Definition~\ref{retractive-def}, but not by Definition~1 of~\cite{bertsekas_tseng_2006}. 

We now state some general properties of the set $\mathcal{R}(X)$. 

\begin{proposition}\label{prop-retractive-cone}
 
For a set $X$, the set $\mathcal{R}(X)$ of retractive directions of $X$ is a nonempty cone.  
  
\end{proposition}

\begin{proof}
If $X=\emptyset$, then $X_\infty=\mathcal{R}(X)=\{0\}.$
Assume now that $X\ne\emptyset$.
Let $d \in \mathcal{R}(X)$ and let $\lambda \geq 0$ be arbitrary. 
Let $\{x_k\} \subset X$ and $\{t_k\} \subseteq \re$ be such that $t_k \rightarrow +\infty$ and $x_k \cdot t_k^{-1} \rightarrow \lambda d$, as $k \rightarrow \infty$, and let $\rho>0$ be arbitrary. Since $d\in\mathcal{R}(X)$ and $\rho\lambda>0$, there exists $K$ such that $x_k-\rho\lambda d\in X$ for all $k\ge K$.
Therefore, $\lambda d\in\mathcal{R}(X)$.
\end{proof}

The cone $\mathcal{R}(X)$ is not necessarily closed, as seen in the following example.

\begin{example}
    Consider the set $X$ given by the epigraph of the function $f(x)=-\sqrt{x}$ for $x\ge0$, i.e., 
    $X=\{(x,\gamma)\in\re^2\mid -\sqrt{x}\le \gamma\}$. The  asymptotic cone of $X$ is the non-negative orthant, i.e.,
    \[X_\infty=\{(d_1,d_2)\in\re^2\mid d_1\ge 0,\ d_2\ge0\}.\]
    We claim that every direction $d=(d_1,d_2)\in X_\infty$ with $d_1 >0$ and $d_2> 0$ is a retractive direction of $X$. To see this, let $\{(x_k, \gamma_k)\}\subset X$ and $\{t_k\} \subseteq \re$ be sequences such that $t_k\to\infty$ and $(x_kt_k^{-1}, \gamma_kt_k^{-1}) \to d$. Let $\rho>0$ be arbitrary. Then, since $(x_k, \gamma_k) \to d$ and $d_1, d_2 >0$, it follows that $x_k \to +\infty$ and $\gamma_k \to+\infty$. Thus, there is a large enough $K$ such that 
    \[x_k - \rho d_1 >0, \quad \text{ and } \quad \gamma_k - \rho d_2 > 0 \qquad\hbox{ for all $k\ge K$}.\]
    Noting that the positive orthant is contained in the set $X$, we see that $(x_k, \gamma_k)-\rho d\in X$ for all $k\ge K$. Hence, $d$ is a retractive direction of the set $X$. 
    
    Next we show that $\mathcal{R}(X)$ is not closed. Note that $(1,0)\in X_\infty$ and consider a sequence $\{d_k\}\subset X_\infty$, with $d_{k,i}>0$ for $i=1,2$, and for all $k$, such that
    \[\lim_{k\to\infty}d_k= (1,0).\]
    As seen above, each $d_k$ is a retractive direction of $X$. However,  the limit $(1,0)$ is not a retractive direction of $X$. To show this, we consider a sequence $\{\bar x_k\}\subset X$ given by 
    \[\bar x_k=(k,-\sqrt{k})\qquad\hbox{for all $k\ge1$},\]
    and note that $\bar x_k\cdot\| \bar x_k\|^{-1}\to (1,0)$.
    For every $k\ge 1$, we can see that 
    \[\bar x_k-(1,0)=(k-1, -\sqrt{k})\notin X.\]
    Hence, $(1,0)$ is not a retractive direction of $X$, implying that $\mathcal{R}(X)$ is not closed.
    \qed
\end{example}

The following proposition considers the cone of retractive directions of the intersection of finitely many sets.

\begin{proposition}\label{prop-retractive-intersect}
   Let $X = \cap_{i=1}^m X_i $ with closed sets $X_i$, for $m\ge 2$. If $X_\infty=\cap_{i=1}^m (X_i)_\infty$, then
   $\cap_{i=1}^m \mathcal{R}(X_i)\subseteq \mathcal{R}(X).$
\end{proposition}

\begin{proof}
    Suppose $d \in \mathcal{R}(X_i)$ for all $i$. Since $\mathcal{R}(X_i)\subseteq (X_i)_\infty$ for all $i$, it follows that $d\in\cap_{i=1}^m (X_i)_\infty.$ By the assumption that $\cap_{i=1}^m (X_i)_\infty= X_\infty$, we have that $d\in X_\infty$.
Let $\{x_k\}\subset X$ be any sequence converging in direction $d$ and let $\rho>0$ be arbitrary. Then, $\{x_k\}\subset X_i$ for all $i$. Since $d\in\mathcal{R}(X_i)$ for all $i$, it follows that for every $i=1,\ldots,m$, there exists an index $K_i$ such that 
    \[x_k-\rho d\in X_i\qquad\hbox{for all }k\ge K_i.\]
   Let $K=\max_{1\le i\le m}K_i$. Then, it follows that 
    \[x_k-\rho d\in X_i\qquad\hbox{for all $k\ge K$ and for all $i=1,\ldots,m$},\]
    implying that
    $x_k-\rho d\in \cap_{i=1}^m X_i$ for all $k\ge K$.
    Thus, $d \in \mathcal{R}(X)$. 
\end{proof}

When the sets $X_i$ in Proposition~\ref{prop-retractive-intersect} are closed and convex, 
the condition $X_\infty=\cap_{i=1}^m (X_i)_\infty$ is always satisfied, as seen from
Proposition~\ref{prop-set-intersect}.

\subsection{Cone of Retractive Directions of Functions}\label{sec-retractive-dirfun}

In this section, we introduce the concept of a \textit{retractive direction of a function}.

\begin{definition}\label{retractive-func}
    An asymptotic direction $d\in\mathcal{K}(f)$ of a proper function $f$ is said to be a retractive direction of $f$
    if for every $\{x_k\}\subset\dom f$ converging in direction $d$ and for every $\rho>0$, there exists an index  $K$ such that 
    \[f(x_k - \rho d) \leq f(x_k)\qquad\hbox{ for all $k \geq K$}. \]
    The set of directions along which $f$ is retractive is denoted by $\mathcal{R}(f)$. 
\end{definition}

By definition,
one can show that $\mathcal{R}(f)$ is a cone and
$0 \in \mathcal{R}(f)$. Furthermore, for a proper function $f$ and a lower-level set $L_\gamma(f)$, we have that 
\[\mathcal{R}(f)\subseteq\mathcal{K}(f)\quad \text{and} \quad  \mathcal{R}(L_\gamma(f))\subseteq\mathcal{K}(f).\]

 It turns out that, in general, there is no special relationship between the cone $\mathcal{R}(f)$ and the cone $\mathcal{R}(L_\gamma(f))$. 
 The following example illustrates that we can have $\mathcal{R}(L_\gamma(f))\subseteq\mathcal{R}(f)$
 for a nonconvex function.

\begin{example}\label{ex-recconef} 
Let $f(s)=\sqrt{|s|}$. Then, 
for every $\gamma>0$, the lower-level set $L_\gamma(f)$ is nonempty and bounded, implying that 
$\mathcal{R}(L_\gamma(f))=\{0\}.$ However, 
we have that 
$(\epi f)_\infty=\{(d,w)\mid d\in\re, w\ge0\}$, implying that 
$f_\infty(d)=0$ for all $d\in\re$.
Moreover, in this case 
$\mathcal{R}(f)=\re$, and we have 
$\mathcal{R}(L_\gamma(f))=\{0\}\subset\mathcal{R}(f)$
for any $\gamma>0$.
\qed
\end{example}

In the following proposition, we establish some properties of $\mathcal{R}(f)$ for a convex function. Interestingly, in this case
$\mathcal{R}(f)\subseteq \mathcal{R}(L_\gamma(f))$ for a lower-level set $L_\gamma(f)$ (converse inclusion to that of Example~\ref{ex-recconef}).

\begin{proposition}\label{prop-rfcones}
    Let $f$ be a proper closed convex function. Then, we have
    \begin{enumerate}
        \item[(a)] \ $\mathcal{R}(f)\subseteq\{d\mid f_\infty(d)=0\}.$
    \item[(b)] $\mathcal{R}(f)\subseteq\mathcal{R}(L_\gamma(f))$ for any $\gamma\in\mathbb{R}$.
    
    \end{enumerate}
\end{proposition}
\begin{proof}
(a) Let $d\in\mathcal{R}(f)$. Then, for any $\{x_k\}$ converging in direction $d$ and any $\rho>0$, there is a large enough $K$ such that
$f(x_k-\rho d)\le f(x_k)$ for all $k\ge K.$
Let $y_k=x_k-\rho d$. Then,  $\{y_k\}$  also converges in direction $d$, and the preceding relation implies that 
\[f(y_k)\le f(y_k+\rho d)\qquad\hbox{for all } k\ge K.\]
By Proposition~2.5.2 of~\cite{auslender_teboulle_2002} for a proper closed convex function $f$, we have that 
\[f_\infty(d)=\sup_{x\in\dom f}\{f(x+d)-f(x)\}.\]
Therefore, it follows that
\[f_\infty(\rho d)\ge \sup_{k\ge K} \{f(y_k+\rho d)-f(y_k)\}\ge 0.\] Since $\rho>0$ and $f_\infty$ is positively homogeneous by Proposition~2.5.1 of~\cite{auslender_teboulle_2002}, it follows that 
$f_\infty (d)\ge 0,$
which combined with the fact that $\mathcal{R}(f)\subseteq \mathcal{K}(f)$, implies that $f_\infty(d)=0$.

\noindent
(b) Let $L_\gamma(f)$ be nonempty. 
 
To arrive at a contradiction,
assume that there is $d\in\mathcal{R}(f)$ such that $d\notin \mathcal{R}(L_\gamma(f))$.
Since $\mathcal{R}(f)\subseteq\mathcal{K}(f)$, it follows that $d\in \mathcal{K}(f)$.
By Proposition~\ref{prop-level-set}, we have that $(L_\gamma(f))_\infty=\mathcal{K}(f)$,
implying that $d\in (L_\gamma(f))_\infty$.
Since $d\in (L_\gamma(f))_\infty$ but it is a non-retractive direction of the set $L_\gamma(f)$, there is a sequence $\{x_k\}\subset L_\gamma(f)$ converging in direction $d$ 
and some $\bar \rho>0$ such that 
\[f(x_k-\bar \rho d)\ge \gamma\qquad\hbox{for infinitely many indices $k$.}\]
Without loss of generality, we may assume that 
$f(x_k-\bar \rho d)\ge \gamma$ for all $k$, 
for otherwise we would choose a suitable subsequence of $\{x_k\}$. Since $\{x_k\}\subset L_\gamma(f) $, it follows that 
\[f(x_k-\bar \rho d)\ge \gamma\ge f(x_k)\qquad\hbox{for all } k.\] 
Thus, we have that $\{x_k\}$ converges in direction $d$ and the direction $d$ belongs to  $\mathcal{K}(f)$, but  $f(x_k-\bar \rho d)\ge f(x_k)$ for all $k$. Hence, the direction $d$ is not retractive for the function $f$,
i.e., $d\not \in \mathcal{R}(f)$, which is a contradiction.
\end{proof}

The inclusion in Proposition~\ref{prop-rfcones}(b) can be strict,
as seen in the following example.
\begin{example}
Let $f(x)=\langle c,x\rangle$ for some $c\in\re^n$, $c\ne0$. Then, for any $\gamma\in\re$, the lower-level set $L_\gamma(f)$ is nonempty and polyhedral, so by Example~\ref{ex-polyhedral-retract} we have that
$\mathcal{R}(L_\gamma(f))=(L_\gamma(f))_\infty$.
Thus, $\mathcal{R}(L_\gamma(f))=\{d\in\re^n\mid \langle c,d\rangle\le 0\}$.

 It can be seen that $f_\infty(d)=\langle c,d\rangle$. By Proposition~\ref{prop-rfcones}, if a direction $d$ is retractive, then  $\langle c,d\rangle =0$.

    Thus, 
    $\mathcal{R}(f)=\{d\mid \langle c,d\rangle=0\}\subset \mathcal{R}(L_\gamma(f))$ for all $\gamma\in\re$.
\qed    
\end{example}

Given two distinct lower-level sets of a function, there is no particular inclusion relation for the cones of their retractive directions, even for a convex function. 
The following examples illustrate that either inclusion between the cones of retractive directions of two lower-level sets is possible.

\begin{example}\label{ex-levelsets-ret}
Consider the function $f(x_1,x_2)=\|x\|-x_1$, for $x_1\ge 0$ and $x_2\in\re$, and level sets
$L_0(f)$ and $L_\gamma (f)$ with $\gamma>0$.
For the set $L_0(f)$ we have
\[L_0(f)=\{(x_1,0)\mid x_1\ge0\},\]
which is a polyhedral set. Thus, $\mathcal{R}(L_0(f))=L_0(f)$.

The set $L_\gamma (f)$ with $\gamma>0$ is given by
$L_\gamma (f)=\{x\in\re^2\mid x_2^2\le 2 x_1\gamma +\gamma^2\}.$
The asymptotic cone of $L_\gamma (f)$ is given by $\{(d_1,0)\in\re^2\mid d_1\ge0\}$ but the retractive direction is only the zero vector,
i.e., $\mathcal{R}(L_\gamma (f))=\{0\}$.
Thus, we have $L_0(f)\subset L_\gamma (f)$, while
$$\mathcal{R}(L_\gamma (f))=\{0\}\subset \mathcal{R}(L_0(f)).$$
        \qed
\end{example}

\begin{example}
Consider the function $f(x_1,x_2)=e^{-\sqrt{x_1 x_2}}$ for $x_1\ge 0$, $x_2\ge0$. 
The lower-level set $L_\gamma (f)$ is nonempty for all $\gamma>0$.
Consider 
$L_1(f)$ and $L_\gamma (f)$ with $\gamma>0$.
For the set $L_1(f)$ we have $L_1(f)=\mathbb{R}^2_+,$
which is a polyhedral set, so we have $\mathcal{R}(L_1(f))=\mathbb{R}^2_+$.
For the set $L_\gamma (f)$ with $\gamma>0$, we have
\[L_\gamma (f)=\{x\in\re^2_+\mid x_1x_2\ge (\ln\gamma)^2\}.\]
The asymptotic cone of $L_\gamma (f)$ is 
$\mathbb{R}^2$  and the cone of retractive directions is 
\[\mathcal{R}(L_\gamma (f))=\{(0,0)\}\cup\{(d_1,d_2)\mid d_1>0, d_2>0\}.\]
Thus, for $\gamma < 1$, we have $L_\gamma(f)\subset L_1 (f)$, while
$$\mathcal{R}(L_\gamma (f))\subset \mathcal{R}(L_1(f)).$$
\qed
\end{example}

Finally, we consider polynomial functions. 
    A polynomial $h: \ren \rightarrow \re$ of order $p$ has the following representation~\cite{Belousov_Klatte_2002}: 
\begin{equation}\label{poly-rep}
    h(x) = \sum_{i=0}^p \phi_i(x) \qquad\hbox{for all } x\in\re^n,
\end{equation}
where each $\phi_i: \ren \rightarrow \re$ is the $i$th order polynomial and $\phi_0$ is a constant. For every $i=0,\ldots,p$,
the polynomial $\phi_i$ has the property that $\phi_i(tx) = t^i\phi_i(x)$ for all $t \in \re$. Asymptotic behavior of polynomials typically depends on their leading order terms. Given a polynomial $h$ of order $p$ and given $x \in \ren$, we let $\mu(x)$ denote the maximal order $i\in\{1,\ldots,p\}$ such that $\phi_i(x)\ne 0$, i.e.,
\[\mu(x)=\max\{ i\mid \phi_i(x)\ne0, i=1,\ldots,p\}.\] 
The following lemma provides a closed form expression for the asymptotic function of a polynomial. 
\begin{lemma}\label{lem-polynomial}

The asymptotic function of a polynomial $h$ of order $p$ is given by 
    \begin{equation*}
        h_\infty(d) = 
        \left\{
        \begin{array}{ll}
            - \infty, &\quad \mu(d)\ge 2\hbox{ and }
            \phi_{\mu(d)}(d) < 0,\\ 
            \phi_1(d), &\quad \mu(d) = 1, \\ 
            + \infty, &\quad \mu(d)\ge 2\hbox{ and } \phi_{\mu(d)}(d) > 0.
        \end{array}
        \right.
    \end{equation*}
    
\end{lemma}

\begin{proof}
    Using the relation $\displaystyle h(x) = \sum_{i=0}^p \phi_i(x)$, for all $x \in \ren$, and the alternative characterization of the asymptotic function, as given in Theorem~\ref{thm-asfun-rep}, we have
    \begin{align*}
            h_\infty(d) = \liminf_{\substack{t \rightarrow \infty \\ d' \rightarrow d}} \frac{h(td')}{t} 
                        &= \liminf_{\substack{t \rightarrow \infty \\ d' \rightarrow d}} \left(\sum_{i=0}^p t^{-1} \phi_i(td')\right) \\ 
                        &= \phi_1(d) + \liminf_{t \rightarrow \infty} \left( \sum_{i=2}^p t^{i-1} \phi_i(d)\right),
        \end{align*} 
        where the last equality follows from $\phi_0/t\to0$, as $t\to\infty$, and the fact that each $\phi_i$ is a continuous function.
    When $\displaystyle \mu(d)\ge 2$,  if $\phi_{\mu(d)}(d) > 0$, then $h_\infty(d)=+\infty$, while  if $\phi_{\mu(d)}(d) < 0$, then $h_\infty(d) = -\infty$. 
    When $\mu(d)=1$, 
   then we are left with $\phi_1(d)$. 
\end{proof}

The cone of retractive directions for a convex polynomial is characterized in the following lemma. It uses a notion of {\it the constancy space of a proper convex function}, defined by
$\mathcal{C}(h)$ is the constancy space of $h$ given by 
\begin{equation}\label{eq-constsp}
\mathcal{C}(h)=\{d\in\re^n\mid f_\infty(d)=f_\infty(-d)=0\}.\end{equation}

\begin{lemma}\label{lem-convpol}
Let $h$ be a convex polynomial of order $p\ge1$. 
Then, we have
\[\mathcal{C}(h)=\{d\mid h_\infty(d)=0\}
\quad\hbox{and} \quad
\mathcal{R}(h)=\{d\mid h_\infty(d)=0\}.\]
\end{lemma}
\begin{proof}
By Lemma~\ref{lem-polynomial} we have that
$$\{d\mid h_\infty(d)=0\}=\{d\mid \phi_1(d)=0, \phi_2(d)=0,\ldots,\phi_p(d)=0\}.$$
Let $d$ be such that $h_\infty(d)=0$. Then,
for every $i=1,\ldots, p$, we have
$\phi_i(-d)=(-1)^i\phi_i(d)$, implying that $h_\infty(-d)=0$. Hence, by Theorem~2.5.3 in~\cite{auslender_teboulle_2002} it follows that \[h(x+t d)=h(x)\qquad\hbox{for all }x\in\dom f.\]
Therefore, the direction $d$ lies in the constancy space $\mathcal{C}(h)$, implying that
\begin{equation}\label{eq-equalsets}
    \{d\mid h_\infty(d)=0\}\subseteq \mathcal{C}(h)
    \subseteq\mathcal{R}(h),
\end{equation}
where the last inclusion in the preceding relation holds since every direction $d\in\mathcal{C}(h)$ is retractive. 
By Proposition~\ref{prop-rfcones}(a), we have that
$\mathcal{R}(h)\subseteq\{d\mid h_\infty(d)=0\}$,
which implies that the equality holds throughout in~\eqref{eq-equalsets}.
\end{proof}

\section{Main Results}\label{sec-main}

In this section, we 
consider the following minimization problem
\begin{equation}\label{opt_problem}
\inf_{x \in X} f(x), \tag{$P$}
\end{equation} 
where $X\subseteq\ren$ is a nonempty closed set and $f$ is a proper closed function.
We let $f^*$ denote the optimal value of the problem, and 
$X^*$ denote the set of optimal solutions.
We provide our sufficient conditions for the existence of solution results in Subsection~\ref{ssec-existence} and relate them to some of the existing  conditions in Subsection~\ref{ssec-related}.

\subsection{Sufficient Conditions for Existence of Solutions}\label{ssec-existence}

The first result relies on the condition that the
asymptotic cone $X_\infty$ of the set 
and the asymptotic cone $\mathcal{K}(f)$ of the function have no nonzero vector in common.

\begin{theorem}\label{thm-compact-case}
Let $X$ be a closed set and $f$ be a proper closed function with $X\cap\dom f\ne \emptyset$.
Assume that $X_\infty\cap {\cal K}(f)=\{0\}.$
Then, the problem \eqref{opt_problem} has a finite optimal value, and  nonempty and compact solution set $X^*$.
\end{theorem}
\begin{proof}
Since $X\cap\dom f\ne \emptyset$, there exists
a point $x_0\in X\cap \dom f$ with a finite value $f(x_0)$. Therefore, for $\gamma=f(x_0)$, the lower-level set $L_\gamma(f)$ is nonempty.
By Proposition~\ref{prop-level-set}, we have that
$(L_\gamma(f))_\infty\subseteq \mathcal{K}(f)$, thus implying that
\[(X\cap L_\gamma(f))_\infty\subseteq X_\infty\cap (L_\gamma(f))_\infty\subseteq X_\infty\cap\mathcal{K}(f)=\{0\},\]
where the first inclusion follows from Proposition~\ref{prop-set-intersect}.
Thus, $(X\cap L_\gamma(f))_\infty=\{0\}$ and the set $X\cap L_\gamma(f)$ is bounded by Proposition~2.1.2 of~\cite{auslender_teboulle_2002}, and hence compact since $X$ and $f$ are closed. 
Therefore, the problem $\inf_{X\cap L_\gamma(f)}f(x)$ has a finite optimal value and a solution exists by the Weierstrass Theorem. Since the problem $\inf_{x\in X\cap L_\gamma(f)} f(x)$ is equivalent to the problem~\eqref{opt_problem},  it follows that $f^*$ is finite and attained. The compactness of $X^*$ follows by noting that $X^*=X\cap L_{f^*}(f)$ is nonempty, closed, and bounded due to $(X\cap L_{f^*}(f))_\infty=\{0\}.$
\end{proof}

For the condition on the asymptotic cones of $X$ and $f$ in Theorem~\ref{thm-compact-case}, we note that 
\[X_\infty\cap\mathcal{K}(f)=\{0\}
\iff f_\infty(0)=0\hbox{ and } f_\infty(d) > 0 \hbox{ for all nonzero }d\in X_\infty.\]
Theorem~\ref{thm-compact-case} generalizes Theorem~3.1 of~\cite{Hieu_2021} 
for minimizing a polynomial function, which requires that $f$ is bounded below on $X$.
It is also more general than Theorem~3 of~\cite{fakhar2023noncoercive},
when restricted to a finite dimensional norm space, which requires that $X$ is convex and $f$ is bounded below on $X$. However, Theorem~3 of~\cite{fakhar2023noncoercive} applies to problems in infinite-dimensional spaces and, by virtue of this, it is more general than our Theorem~\ref{thm-compact-case}.

We make a final remark that the condition of~Theorem~\ref{thm-compact-case} implies that 
\[(f+\delta_X)_\infty(d)>0\qquad \hbox{for all nonzero $d\in\re^n$},\]
where $\delta_X$ is the indicator function of the set $X$.

Thus, the function $f+\delta_X$ satisfies the conditions of Theorem~3.4.1 in~\cite{auslender_teboulle_2002}, which implies that the problem~\eqref{opt_problem} has a solution, but it does not imply that the solution set is bounded.

We now provide sufficient conditions for the existence result using the asymptotically bounded decay property with respect to a coercive function $g$.

\begin{theorem}\label{main-result-2g}
    Let the set $X \subseteq \ren$ be closed, and let the function \propf be proper closed and and $X \cap \dom f \neq \emptyset$.
    Suppose that $f$ exhibits asymptotically bounded decay with respect to a coercive function $g$ on $X$. Assume that 
    $$X_\infty \cap \mathcal{K}(f) \subseteq \mathcal{R}(X) \cap \mathcal{R}(f).$$ 
    Let $g$ have Lipschitz continuous gradients on $X$ and, for any sequence $\{x_k\} \subseteq X$ converging in any nonzero direction $d \in \mathcal{R}(X) \cap \mathcal{R}(f)$, let the following relation hold
    \[\liminf_{k\to\infty} \frac{\langle \nabla g(x_k),d\rangle}{t_k}>0\qquad\hbox{where } \lim_{k\to\infty} \frac{x_k}{t_k}=d.\]  
    Then, the problem (\ref{opt_problem}) has a finite optimal value $f^*$ and a nonempty solution set $X^*$. 
    \end{theorem}

    \begin{proof}
    By our assumption that
    $f$ exhibits asymptotically bounded decay with respect to $g$ on $X$,
    we have by Definition \ref{def-asymp-decay} that
    \begin{equation}
    \label{eq-decay-new}
    c=\liminf_{\substack{\|x\| \rightarrow \infty \\ x \in X}} \frac{f(x)}{g(x)} >-\infty.
 \end{equation} 
 If $c=+\infty$, then $f$ is coercive since $g$ is such. Thus, the result follows by  Proposition~\ref{prop-exist-basic}.
 
   Next,  we consider the case $c\in\re$. For this,
  we consider the regularized problem
        \begin{equation}\label{eq-reg-g}
            \inf_{x \in X} \{f(x) + r g^2(x)\}
            \qquad\hbox{for any $r>0$}.
        \end{equation}
  \noindent
        (\textit{Step 1: The regularized problem has a solution.}) 
        Since $g$ is coercive, there exists an $R_g>0$ such that 
        \[g(x)>0 \qquad\hbox{for all $x\in X$ with $\|x\|\ge R_g$.}\]
        By \eqref{eq-decay-new}, for any $\epsilon > 0$, there exists large enough $R_\epsilon> R_g> 0$ such that  
        $$ f(x) \geq (c-\epsilon)g(x)\qquad\hbox{for all $x\in X$ with $\|x\|\ge R_\epsilon$}.$$ 
        Thus, for all $x\in X$ with $\|x\|\ge R_\epsilon$,
        $$f(x) + r g^2(x)\geq (c - \epsilon)g(x) + r g^2(x),$$ implying that
        $$ \liminf_{\substack{\|x\| \rightarrow \infty \\ x \in X}} \{f(x) + r g^2(x)\} 
        \geq \liminf_{\substack{\|x\| \rightarrow \infty \\ x \in X}} g(x)\left(rg(x)+c-\epsilon)\right)=
        +\infty, $$
        where the equality is due to $r>0$ and
        the coercivity of $g$ on $X$. Hence, the regularized objective function itself is also coercive, and by Proposition~\ref{prop-exist-basic} the problem $\inf_{x\in X} \{f(x) + r g^2(x)\}$ has a solution for any $r>0$. 
        
        \noindent
        (\textit{Step 2: A sequence of solutions to regularized problems is bounded.})
        Now, consider a sequence of positive scalars $\{r_k\}$ such that $r_k \rightarrow 0$ as $k \to \infty$.
        For each $k$, let $x_k^* \in X$ be a solution to the regularized problem in~\eqref{eq-reg-g} with $r=r_k$.
        Towards a contradiction, assume that the sequence $\{x_k^*\} \subseteq X$ is unbounded. 
        
         Without loss of generality, let $x_k^*\ne 0$ for all $k$, and consider the sequence $\{x_k^*/\|x_k^*\|\}$. This sequence is bounded and, hence, must have a convergent subsequence. Again without loss of generality, let $\{x_k^*/\|x_k^*\|\} \rightarrow d$ as $k \to \infty$. Thus, $d \in X_\infty$. 
         Since $x_k^*\in X$ is a solution of the regularized problem, we have
        \begin{equation}\label{eq-fandg} 
            f(x_k^*) \leq f(x_k^*) + r_k g^2(x_k^*) \leq f(x) + r_k  g^2(x) \quad \text{ for all $x\in X$ and all }  k. 
        \end{equation}
        Thus, it follows that for an arbitrary fixed $x_0\in X$,
         \begin{equation*} 
            \liminf_{k\to\infty}\frac{f(x_k^*)}{\|x_k\|} \leq \lim_{k\to\infty}\frac{f(x_0)+ r_k  g^2(x_0)}{\|x_k\|} =0, 
        \end{equation*}
        where we use the fact $r_k\to0$.
        Employing Theorem~\ref{thm-asfun-rep} which gives the explicit form for the asymptotic function $f_\infty$, we can see that
    $$f_\infty(d) \leq \liminf_{k\rightarrow \infty} \frac{f(\|x_k^*\|(\|x_k^*\|^{-1}x_k^*))}{\|x_k^*\|} = \liminf_{k\rightarrow\infty} \frac{f(x_k^*)}{\|x_k^*\|} \leq 0. $$ 
      Hence $d\in X_\infty \cap \mathcal{K}(f)$. By the assumption that $X_\infty \cap \mathcal{K}(f)\subseteq \mathcal{R}(X) \cap \mathcal{R}(f)$, it follows that $d \in \mathcal{R}(X) \cap \mathcal{R}(f)$. Since $d$ is a retractive direction of the set $X$ and the function $f$, for the sequence $\{x_k^*/\|x_k^*\|\}$ converging in the direction $d$ and any for $\rho > 0$, there exists a sufficiently large index $K$ so that 
      $$ x_k^* - \rho  d \in X \quad \text{ and } \quad f(x_k^* - \rho d) \leq f(x_k)\qquad\hbox{for all $k \geq K$}. $$
      Thus, for any $k \geq K$, 
        \begin{equation}\label{eq-greg}
            f(x_k^*) + r_k g^2(x_k^*)\leq f(x_k^* - \rho d) + r_kg^2(x_k^* - \rho d) \leq f(x_k^*) + r_k g^2(x_k^*-\rho d)^2,
        \end{equation}
        where the first inequality follows from the optimality of the point $x_k^*$ for the regularized problem and the fact that $x_k-\rho d\in X$, while 
    the second inequality follows by $f(x_k^*-\rho d)\le f(x_k^*)$.
         Since $r_k>0$ for all $k$, relation~\eqref{eq-greg} implies that 
         \begin{align}\label{eq-similar}
         g^2(x_k^*) \leq g^2(x_k^* - \rho d)\qquad\hbox{for all $k \geq K$}.
         \end{align}
         Assume that $K$ is large enough so that $g(x_k)>0$ for all $k\ge K$.
        Then, relation~\eqref{eq-similar} implies that 
        \begin{equation}\label{contra-implication}
           0< g(x_k^*) \leq g(x_k^* - \rho d) \le g(x_k^*)-\rho\langle \nabla g(x_k^*),d\rangle +\frac{\rho^2 L}{2}\qquad \text{ for all } k \ge K.
        \end{equation} 
        where in the last inequality we use the Lipschitz continuity of $\nabla g$ and  $\|d\|=1$.
        Since $\rho>0$, by relation~\eqref{contra-implication} it follows that 
        \[\langle \nabla g(x_k^*),d\rangle \le \frac{\rho L}{2}\qquad \text{ for all } k \ge K.\]
        Hence, 
        \[\liminf_{k\to\infty}\frac{\langle \nabla g(x_k^*),d\rangle}{\|x_k^*\|} 
        \le \liminf_{k\to\infty}\frac{\rho L}{2\|x_k^*\|}=0,\]
        which contradicts the assumption on $\nabla g$ that $\liminf_{k\to\infty}\frac{\langle \nabla g(x_k^*),d\rangle}{\|x_k^*\|} >0$.
        Thus, the sequence $\{x_k^*\}$ must be bounded. 

        \noindent
        (\textit{Step 3: Any accumulation point of the sequence $\{x_k^*\}$ is a solution to \eqref{opt_problem}})
        Since the sequence $\{x_k^*\}$ is bounded, it must have an accumulation point. Without loss of generality, let $\{x_k^*\} \rightarrow x^*$ as $k \to \infty$. Taking the limit inferior in~\eqref{eq-fandg} yields
        \begin{equation*} 
            \liminf_{k \to\infty}f(x_k^*) \leq  \liminf_{k\to\infty}\{f(x) + r_k  g^2(x) \}\quad \text{ for all $x\in X$.} 
        \end{equation*}
        Since $x_k^* \to x^*$ and $r_k \to 0$,  by the closedness of $f$  it follows that 
        $$ f(x^*) \leq f(x) \quad \text{ for all } x \in X.$$
        Since $\{x_k^*\} \subseteq X$ and $X$ closed,  $x^* \in X$ is an optimal solution to the problem \eqref{opt_problem}. 
    \end{proof}

The proof of Theorem~\ref{main-result-2g}
is motivated by the insights gained from 
proof of
Theorem~3.4.1 in~\cite{auslender_teboulle_2002}, which is in turn motivated by~\cite{baiocchi1998}, and it has ideas similar to work~\cite{penot2006} on the existence of solutions for problems (in a Banach space).

The following result is a consequence of 
   
    Theorem~\ref{main-result-2g} for a polyhedral set $X$.
    
    \begin{corollary}\label{cor-polyhedral-set}
    Let $X$ be a polyhedral set and $f$ be a proper closed function with $X\cap\dom f\ne\emptyset$.
     Under assumptions of Theorem~\ref{main-result-2g}, if 
     $X_\infty\cap\mathcal{K}(f)\subseteq X_\infty\cap \mathcal{R}(f),$
     then the problem~\eqref{opt_problem} has a solution.
    \end{corollary}
    \begin{proof}
        The result follows by Theorem~\ref{main-result-2g} since a  polyhedral set is retractive i.e., $\mathcal{R} (X)=X_\infty$ (see Example~\ref{ex-polyhedral-retract}).
    \end{proof}

\begin{remark}
\label{rem-conds}
The condition $X_\infty\cap\mathcal{K}(f)\subseteq \mathcal{R}(X)\cap \mathcal{R}(f)$ of 
Theorem~\ref{main-result-2g} implies that we must have 
\[X_\infty \cap\mathcal{K}(f)= \mathcal{R}(X) \cap \mathcal{R}(f),\]
since 
$\mathcal{R}(X)\cap \mathcal{R}(f)\subseteq X_\infty\cap\mathcal{K}(f)$ always holds by the definition of the retractive cones.
\end{remark} 

We next provide a result for the special case when the conditions of Theorem~\ref{main-result-2g} hold with $g(x)=\|x\|^p$ for some $p\ge0$. In this case, we can have $p=0$ so that $g$ need not be coercive, and the Lipschitz continuity of $\nabla g$
is not needed.

\begin{theorem}\label{main-result}
    Let the set $X \subseteq \ren$ be closed, and let the function \propf be proper and closed with $X \cap \dom f \neq \emptyset$. Assume that 
    $f$ exhibits asymptotically bounded decay with respect to $g(x) = \|x\|^p$ on $X$ for some $p \geq 0$, and assume that 
    \[X_\infty \cap \mathcal{K}(f) \subseteq \mathcal{R}(X) \cap \mathcal{R}(f).\] Then, the problem (\ref{opt_problem}) has a finite optimal value $f^*$ and a nonempty solution set $X^*$.
\end{theorem}

\begin{proof}
The proof is along the same lines as that of Theorem~\ref{main-result-2g}, which we follow and point to simplifications due to the choice $g(x)=\|x\|^p$.
By our assumption that
    $f$ exhibits asymptotically bounded decay with respect to $g(x) = \|x\|^p$ on $X$ for some $p \geq 0$,
    we have that
    \begin{equation}
    \label{eq-decay}
    c=\liminf_{\substack{\|x\| \rightarrow \infty \\ x \in X}} \frac{f(x)}{\|x\|^p} >-\infty.
 \end{equation} 
 If $c=+\infty$, then $f$ is coercive on $X$ so by Proposition~\ref{prop-exist-basic}, the problem~(\ref{opt_problem}) has an optimal solution. So, consider the case $c\in\re$.

\noindent    
(\textit{Step 1: The regularized problem has a solution.})
Let $\epsilon > 0$ be arbitrarily small. By the asymptotically bounded decay of $f$ on $X$ in~\eqref{eq-decay}, there exists $R>0$ large enough so that
    $$f(x) \geq (c-\epsilon)\|x\|^p \qquad\hbox{for all $x\in X$ with $\|x\|\ge R$}.$$
    Let $r>0$ be arbitrary and consider the function $f(x)+r\|x\|^{p+1}$. We have
    $$f(x) + r\|x\|^{p+1}  \geq r \|x\|^{p+1} + (c-\epsilon)\|x\|^p
    \qquad\hbox{for all $x\in X$ with $\|x\|\ge R$},$$ 
    implying that 
    $$\liminf_{\substack{\|x\| \rightarrow \infty \\ x \in X}}  \left\{f(x) +r\|x\|^{p+1} \right\} 
    \geq \liminf_{\substack{\|x\| \rightarrow \infty \\ x \in X}}  \left\{\|x\|^p\left(r\|x\| +c-\epsilon\right) \right\} 
    =+ \infty. $$ Thus, for any $\rho>0$, the function $f(x) + r\|x\|^{p+1}$ is coercive, so by Proposition~\ref{prop-exist-basic},  
    the regularized problem 
    $\inf_{x \in X} \{f(x) + r\|x\|^{p+1}\}$
    has a solution for every $r>0$. 

\noindent
    (\textit{Step 2: A sequence of solutions to regularized problems is bounded.})
    This part is the same as Step 2 of the proof of Theorem~\ref{main-result-2g}, where $g^2(x)$ is replaced with $\|x\|^{p+1}$ up to and including relation~\eqref{eq-similar},
    which in this case reduces to 
    $$ \|x_k^*\|^{p+1} \le \|x_k^* - \rho d\|^{p+1}\qquad\hbox{for $k \geq K$}. $$ 
    Therefore,
    \begin{equation*}\label{contra-ineq}
        \|x_k^*\|^2 \leq \|x_k^* - \rho d\|^2
        \qquad\hbox{for all $k \geq K$}. 
    \end{equation*} 
    Since $\|d\|=1$, the preceding inequality implies that $ 2\langle x_k^*, d \rangle \leq \rho$ for all $k \geq K$. Therefore, it follows that
    $$ \lim_{k\rightarrow \infty} \frac{\langle x_k^*, d \rangle}{\|x_k^*\|} \leq \lim_{k\rightarrow \infty} \frac{\rho}{2\|x_k^*\|} = 0.$$
    On the other hand, since $\lim_{k\to\infty}x_k^*\cdot\|x_k^*\|^{-1}=d$ and $\|d\|=1$, we obtain that
    
    $1\le 0$,
    which is a contradiction.
    Hence, it must be that $\{x_k^*\}$ is bounded. 
    
\noindent
    (\textit{Step 3: Any accumulation point of sequence $\{x_k^*\}$ is a solution to~\eqref{opt_problem}.}) The proof is the same as Step 3 of the proof of Theorem~\ref{main-result-2g}, where $g^2(x)$ is replaced with $\|x\|^{p+1}$.
    \end{proof}

    Next, we provide an example showing that, if the condition $X_\infty\cap \mathcal{K}(F)\subseteq \mathcal{R}(X)\cap\mathcal{R}(f)$ is violated, the problem may not have a finite optimal value; hence no solution.

\begin{example}\label{ex-noval}
Consider the problem of minimizing a convex scalar function $f(x)=-\sqrt{x}$ over its domain $X=\{x\mid x\ge0\}$.
The optimal value is $f^*=-\infty$ and there is no solution.
Since the function is convex, by Example~\ref{ex-convex}, it exhibits asymptotically bounded decay with respect to $g(x)=|x|$ on $X$. 
The set $X$ is a closed convex cone, and we have $X_\infty=X$.
The asymptotic cone of $f$ coincides with $X$ i.e., $\mathcal{K}(f)=X$, implying that
$X_\infty\cap \mathcal{K}(f)=X.$
The cone of retractive directions of $X$ coincides with $X$, since $X$ is a polyhedral set (see Example~\ref{ex-polyhedral-retract}).
The cone $\mathcal{R}(f)$ of retractive directions of $f$ contains only the zero vector. To see this note that for $x_k=\lambda k$, with $\lambda>0$, and $t_k=k$ for all $k\ge 1$, we have that $\{x_k\}$ converges in the direction $d=\lambda$. However, for any $\rho\in(0,1)$ and any $k\ge 1$,
\[f(x_k-\rho \lambda)=-\sqrt{k-\rho\lambda}>-\sqrt{k}=f(x_k).\]
Hence, we have $\mathcal{R}(f)=\{0\}$, and 
$\mathcal{R}(X)\cap\mathcal{R}(f)=\{0\},$
thus implying that the condition $X_\infty\cap \mathcal{K}(f)\subseteq \mathcal{R}(X)\cap\mathcal{R}(f)$ of Theorem~\ref{main-result} is violated.
\end{example}

The following example shows that when the condition $X_\infty \cap \mathcal{K}(f)\subseteq \mathcal{R}(X) \cap \mathcal{R}(f)$ is violated, the problem  can have a finite optimal value but not a solution.

\begin{example}
Consider the problem of minimizing a convex scalar function $f(x)=e^x$ over its domain $X=\re$.
The optimal value is $f^*=0$ and there is no solution.
The function is convex, so by Example~\ref{ex-convex}, it exhibits asymptotically bounded decay with respect to $g(x)=|x|$ on $X$. 
We have $X_\infty=X$ and $\mathcal{R}(X)=X$.
The asymptotic cone of $f$ coincides is given by
\[\mathcal{K}(f)=\{x\mid x\le 0\},\]
while
the cone $\mathcal{R}(f)$ of retractive directions of $f$ contains only the zero vector. To see this, observe that for $x_k=-\lambda k$, with $\lambda>0$, and $t_k=k$ for all $k\ge 1$, we have that $\{x_k\}$ converges in the direction $d=-\lambda$. However, for any $\rho\in(0,1)$ and any $k\ge 1$,
\[f(x_k+\rho \lambda)=e^{k+\rho\lambda}>e^{k}=f(x_k).\]
Hence, we have $\mathcal{R}(f)=\{0\}$ and $\mathcal{R}(X)\cap\mathcal{R}(f)=\{0\},$
implying that the condition $X_\infty\cap \mathcal{K}(f)\subseteq \mathcal{R}(X)\cap\mathcal{R}(f)$ of Theorem~\ref{main-result} is violated.
\end{example}

As the preceding examples show, when the condition on the retractive cones of Theorem~\ref{main-result} are violated, the problem~\eqref{opt_problem} may not have finite optimal value $f^*$, or it may have finite $f^*$ but not a solution. This leads us to conjecture that the sufficient conditions of Theorem~\ref{main-result}, and Theorems~\ref{thm-compact-case}--\ref{main-result-2g} might also be necessary. It turns out that it is indeed the case, as 
seen in the next section.

\subsection{Conditions of Theorems~\ref{thm-compact-case}--\ref{main-result} are also Necessary}\label{ssec-related}
To show that the conditions of Theorems~\ref{thm-compact-case}--\ref{main-result}  are also necessary for the existence of solutions to problem~\eqref{opt_problem}, we use Corollary~4.1 in~\cite{penot2006}, which provides necessary and sufficient conditions for solution existence to problems in a Banach space.
We restate that corollary for the case of $\mathbb{R}^n$.

\begin{proposition}[Corollary~4.1 in~\cite{penot2006}]\label{prop-penot}
Let the set $X \subseteq \ren$ be closed, and let the function \propf be proper and closed.
    Suppose that $h$ is a coercive function and that $X \cap \dom f \cap \dom h \neq \emptyset$. 
    Then, the following two conditions are necessary and sufficient in order that the function $f$ attains its infimum on $X$:
    \begin{itemize}
     \item[a)]  For any (nonzero) $d\in X_\infty$, for any sequence $\{x_k\}\subset X$,
     such that $f(x_k)\to \inf_{x\in X}f(x)$ and $x_k/\|x_k\|\to d$, for all $k$ large enough there exists $y_k\in X$ such that
     \[f(y_k)\le f(x_k),\qquad h(y_k) < h(x_k);\]
\item[b)] The function $f$ is bounded below on the set $X$.
    \end{itemize}
\end{proposition}
We note that when $X_\infty=\{0\}$ then the condition (a) is vacuous and, thus, satisfied by default.
Using Proposition~\ref{prop-penot}, we obtain the following result.

\begin{theorem}\label{thm-necessity}
The conditions of Theorems~\ref{thm-compact-case}--\ref{main-result} are also necessary for the existence of solutions to problem~\eqref{opt_problem}.
\end{theorem}
\begin{proof}
 
To show the necessity, we prove that the conditions of each of Theorems~\ref{thm-compact-case}--\ref{main-result} imply that the conditions of Proposition~\ref{prop-penot} are satisfied. 
Since the conditions of Proposition~\ref{prop-penot} are necessary, the conditions of Theorems~\ref{thm-compact-case}--\ref{main-result} are also necessary.

We note that under the conditions of Theorems~\ref{thm-compact-case}--\ref{main-result}, we have that $f^*=\inf_{x\in X}f(x)$ is finite. Hence, $f$ is bounded below on $X$, and 
the condition (b) of Proposition~\ref{prop-penot} is satisfied.
It remains to show that the condition (a) of Proposition~\ref{prop-penot} is satisfied. For this, we let 
$d\in X_\infty$ and  $\{x_k\}\subset X$ be any sequence
     such that $f(x_k)\to \inf_{x\in X}f(x)$ and $x_k/\|x_k\|\to d$. Then, $d\in X_\infty$ and, since $\{x_k\}\subset X$,
     \[f_\infty(d)\le \liminf_{k\to\infty}\frac{f(x_k)}{\|x_k\|}\le0,\]
     where the last inequality follows from either $\inf_{x\in X}f(x)$ being finite or equal to $-\infty$. Hence, $d\in\mathcal{K}(f)$, and consequently 
     \[d\in X_\infty \cap \mathcal{K}(f) .\]
Now, we consider each of Theorems~\ref{thm-compact-case}--\ref{main-result} separately.\\
\noindent
{\it Case of Theorem~\ref{thm-compact-case}}: The condition $X_\infty \cap \mathcal{K}(f)=\{0\}$
implies that we must have $d=0$. Thus, there is no nonzero vector satisfying the condition (a) of Proposition~\ref{prop-penot}. Hence, the condition is satisfied by default.\\
\noindent
{\it Case of Theorems~\ref{main-result-2g}--\ref{main-result}}:
By the condition
    $X_\infty \cap \mathcal{K}(f) \subseteq \mathcal{R}(X) \cap \mathcal{R}(f)$
    of Theorems~\ref{main-result-2g}--\ref{main-result}, and the definitions of cones of retractive directions for $X$ and $f$, 
    for the sequence $\{x_k\}\subset X$, with $x_k/\|x_k\|\to d$, and for any $\rho>0$, there exists $K$ (depending on $\rho$) such that 
    \[x_k-\rho d\in X,\qquad f(x_k-\rho d)\le f(x_k)\qquad\hbox{for all }k\ge K.\]
    Hence, the condition $f(y_k)\le f(x_k)$ in (a) of Proposition~\ref{prop-penot} is satisfied with $y_k=x_k-\rho d\in X$ (for an arbitrary $\rho>0$.)

We next show that the condition $h(y_k)< h(x_k)$ is also satisfied with $y_k=x_k-\rho d$ for a larger index $K$, for the function $h(x)=\|x\|^2$.

Since $\|d\|=1$,  we have
 \[\|x_k-\rho d\|^2 = \|x_k\|^2 +\rho^2 -2\rho\langle x_k,d\rangle
 =\|x_k\|^2 +\rho\|x_k\|\left(\frac{\rho}{\|x_k\|} -2\left\langle \frac{x_k}{\|x_k\|},d\right\rangle\right).\]
 Since $\|x_k\|\to\infty$ and $x_k/\|x_k\|\to d$, it follows that 
 \[\lim_{t\to\infty} \left(\frac{\rho}{\|x_k\|} -2\left\langle \frac{x_k}{\|x_k\|},d\right\rangle\right)=-2.\]
 Thus, there exists an index $K_1\ge K$ and an $\epsilon\in(0,2)$ so that 
\[\frac{\rho}{\|x_k\|} -2\left\langle \frac{x_k}{\|x_k\|},d\right\rangle\le \epsilon-2 <0
 \qquad \hbox{for all }k\ge K_1.\]
 Hence, for $k\ge K_1$ and $y_k=x_k-\rho d$, we have 
 $\|y_k\|^2 <\|x_k\|^2$ for all $k\ge K_1$. Thus, the condition $h(y_k)< h(x_k)$ in (a)
 of Proposition~\ref{prop-penot} is satisfied with $h(x)=\|x\|^2$.
\end{proof}

By Theorem~\ref{thm-necessity} and Corollary~\ref{cor-polyhedral-set}, the following result holds when $X$ is polyhedral.
\begin{corollary}\label{cor-polyhedral-set-nec}
     The assumptions of Corollary~\ref{cor-polyhedral-set} are necessary for 
     the problem~\eqref{opt_problem}, with a polyhedral set $X$, to have a solution.
    \end{corollary}

\section{Implications for Convex Problems}\label{sec-conv}

In this section, we discuss the relationship between the asymptotic cone and the lineality space~\cite{bertsekas_nedic_ozdaglar_2013} of a closed convex set. 

Then, we apply Theorem~\ref{main-result-2g} to provide necessary and sufficient  conditions for existence of solutions to a convex minimization problem.

For a nonempty closed convex set $X\subseteq\re^n$,
the asymptotic cone $X_\infty$ has a simple
characterization (Proposition 2.1.5 of~\cite{auslender_teboulle_2002})
\[X_\infty=\{d\in\re^n\mid \exists x\in X\hbox{ such that }x+t d\in X,\,\forall t\ge0\}.\]

The asymptotic cone $X_\infty$ is often referred to as a recession cone of $X$~\cite{rockafellar-1970a,bertsekas_nedic_ozdaglar_2013}.
The lineality space of a nonempty closed convex set $X$
is the set
defined as ${\rm Lin}(X)=X_\infty\cap(-X_\infty)$ ~\cite{rockafellar-1970a,bertsekas_nedic_ozdaglar_2013}, and satisfies
\[{\rm Lin}(X)=\{d\in\re^n\mid \exists x\in \re^n\hbox{ such that } x+t d \in X,\,\forall t\in\re\}.\]

For a nonempty closed convex set $X$, by the definition of the cone $\mathcal{R}(X)$ of retractive directions of $X$,
we have
\[{\rm Lin}(X)\subseteq \mathcal{R}(X).\]
The inclusion can be strict. For example, if $X=\{x\in\re^n\mid Ax\le b\}$ for some matrix $A$ and a vector $b$, then ${\rm Lin}(X)=\{d\mid Ad=0\}$ and $\mathcal{R}(X)=X_\infty=\{d\mid Ad\le 0\}$ since $X$ is polyhedral.

The constancy space $\mathcal{C}(f)$ (cf.~\eqref{eq-constsp}) of a proper closed convex function $f$ satisfies the following relations (see Theorem~2.5.3 in~\cite{auslender_teboulle_2002}):
\[\mathcal{C}(f)=\{d\in\re^n\mid \exists x\in \dom f\hbox{ such that }f(x+t d)= f(x),\,\forall t\in\re\},\]
or equivalently
$\mathcal{C}(f)=\{d\in\re^n\mid f(x+t d)= f(x),\,\forall x\in\dom f, t\in\re\}.$
For a proper closed convex function $f$, by the definition of the cone $\mathcal{R}(f)$ of retractive directions of $f$, it follows that 
\[\mathcal{C}(f)\subseteq \mathcal{R}(f).\]

\subsection{Necessary and Sufficient Conditions for Solution Existence}

We next consider a general convex problem of the form
\begin{eqnarray}\label{convex-prob}
\hbox{minimize \ \,} && f(x)\cr 
\hbox{subject to } && g_j(x)\le 0, j=1,\ldots,m,\ x\in C.
\end{eqnarray}
We have the following result.

\begin{theorem}\label{thm-convex}
Let $X$ be a closed convex set, and let $f$ and each $g_j$ be a proper closed   convex function such that 
$C\cap (\cap_{j=1}^m \dom g_j)\cap\dom f\ne\emptyset$. Then, the problem~\eqref{convex-prob} has a finite optimal value and its solution set is nonempty if and only if
\begin{equation}\label{eq-convex-cond}
    C_\infty \cap \left(\cap_{j=1}^m \mathcal{K}(g_j)\right) \cap \mathcal{K}(f)\subseteq {\cal R}(C) \cap\left( \cap_{j=1}^m{\cal R}(L_0(g_j))\right)\cap {\cal R}(f).
\end{equation}
Moreover, if the set $C$ is polyhedral, then 
the problem~\eqref{convex-prob} has a finite optimal value and its solution set 
set is nonempty if and only if  
\begin{equation}\label{eq-convex-cond1}
    C_\infty \cap\left(\cap_{j=1}^m \mathcal{K}(g_j)\right) \cap \mathcal{K}(f) \subseteq C_\infty\cap\left( \cap_{j=1}^m{\cal R}(L_0(g_j))\right)\cap {\cal R}(f).
\end{equation}
\end{theorem}
\begin{proof}
Since the objective function is convex, it exhibits asymptotically bounded decay with respect to $g(x)=\|x\|$ (see Example~\ref{ex-convex}). 
Furthermore, 
let 
\[X=\{x\in C\mid g_j(x)\le 0,\ j=1,\ldots,m\}.\]
The set $X$ is nonempty, closed, and convex, so by Proposition~\ref{prop-set-intersect}, we have 
$$X_\infty=C_\infty\cap\left(\cap_{j=1}^m\{x\in\ren\mid g_j(x)\le 0\}_\infty\right).$$
By Proposition~\ref{prop-level-set}, we have that 
$\{x\mid g_j(x)\le 0\}_\infty=\mathcal{K}(g_j)$ for all $j=1,\ldots,m$,
thus implying that
$$X_\infty=C_\infty\cap\left(\cap_{j=1}^m \mathcal{K}(g_j)\right).$$
Moreover, by Proposition~\ref{prop-retractive-intersect}, we have that 
$${\cal R}(C) \cap\left( \cap_{j=1}^m{\cal R}(L_0(g_j))\right)\subseteq \mathcal{R}\left( C \cap (\cap_{j=1}^m L_0(g_j))\right)=\mathcal{R}(X).$$
The preceding two relations combined with~\eqref{eq-convex-cond} show that the condition
$X_\infty \cap \mathcal{K}(f) \subseteq \mathcal{R}(X) \cap \mathcal{R}(f)$
of Theorem~\ref{main-result} is satisfied, and the sufficiency part follows 
by Theorem~\ref{main-result} and the necessary part by Theorem~\ref{thm-necessity}.

When  $C$ is polyhedral, the result follows by Corollary~\ref{cor-polyhedral-set-nec} and Theorem~\ref{thm-necessity}. 
\end{proof}

Theorem~\ref{thm-convex} is more general than Proposition~6.5.4 in~\cite{bertsekas_nedic_ozdaglar_2013}, which provides only sufficient conditions for the existence of solutions, under more stringent assumptions than those of Theorem~\ref{thm-convex} 
requiring that problem~\eqref{convex-prob} has a finite optimal value and
\[C_\infty \cap \left(\cap_{j=1}^m \mathcal{K}(g_j)\right) \cap \mathcal{K}(f)\subseteq 
{\rm Lin}(C)\cap\left(\cap_{j=1}^m \mathcal{C}(g_j)\right) \cap\mathcal{C}(f).\]
The preceding condition implies that
the condition~\eqref{eq-convex-cond} holds since ${\rm Lin}(C)\subseteq\mathcal{R}(C)$, and the analogous relation holds for the constancy space and the cone of retractive directions for the functions $f$ and $g_j$. 
Similarly, for the case of a polyhedral set $C$, Theorem~\ref{thm-convex} 
is more general than Proposition~6.5.5 in~\cite{bertsekas_nedic_ozdaglar_2013}.

Moreover, when the set $C$ is polyhedral and $g_j\equiv 0$ for all $j$, 
the condition~\eqref{eq-convex-cond1} reduces to
\[C_\infty \cap \mathcal{K}(f) \subseteq C_\infty \cap {\cal R}(f).\]
In this case, Theorem~\ref{thm-convex} is more general than Theorem~27.3 in~\cite {rockafellar-1970a}, which 
provides only sufficient conditions using a more restrictive condition
that $C_\infty \cap \mathcal{K}(f) \subseteq 
\mathcal{C}(f)$.

\section{Implications for Nonconvex Problems}
\label{sec-nonconv}

In this section, we consider the implications of our main results for several types of
nonconvex problems for the cases where the constraint set $X$ is generic and given by nonconvex functional inequalities (Sections~\ref{sec-esnonalg1} and~\ref{sec-algsets}, respectively).

\subsection{Generic Constraint Set}
\label{sec-esnonalg1}

We consider the problem~\eqref{opt_problem} for the case when $f$ is convex while the set $X$ need not be convex, for which we have by Theorem~\ref{main-result} and Theorem~\ref{thm-necessity}.

\begin{theorem}\label{thm-conf-nonconvx}
Let $X$ be a nonempty closed set and $f$ be a proper closed convex function with $X\cap\dom f\ne\emptyset$.
Then, the problem~\eqref{opt_problem} has a finite optimal value and a solution exists
if and only if
\[X_\infty\cap \mathcal{K}(f)\subseteq \mathcal{R}(X)\cap \mathcal{C}(f).\]
\end{theorem}
\begin{proof}
A convex function exhibits asymptotically bounded decay with respect to $g(x)=\|x\|$ (Example~\ref{ex-convex}), and we have $\mathcal{C}(f)\subseteq\mathcal{R}(f)$ when $f$ is proper closed convex function.
The sufficiency part follows by Theorem~\ref{main-result}, while the necessary part follows from  Theorem~\ref{thm-necessity}.
\end{proof}

To the best of our knowledge the result of Theorem~\ref{thm-conf-nonconvx} is new. An existing result that considers convex objective and a nonconvex constraint set is Proposition~12 in~\cite{bertsekas_tseng_2006}, which relies on the stringent assumption that:\\
\noindent (A1) Every nonzero direction $d \in X_\infty$ is retractive and, for all $x \in X$, there exists an $\bar \alpha \geq 0$ such that  $x + \alpha d \in X$ for all  $\alpha \geq \bar \alpha.$

In the following example assumption (A1) fails to hold. Thus, Proposition~12 in~\cite{bertsekas_tseng_2006} cannot be applied to assert the existence of solutions,
while Theorem~\ref{thm-conf-nonconvx} can be used. 
    
    \begin{example}\label{ex-ourthm-works}
        Consider minimizing a proper closed convex function $f$ on the set $X$ given by 
        $X = \{x\in\re^2 \mid x_2 \leq x_1^2 \}$
        (cf.\ Figure \ref{fig:asymp-cone-neq-recc}).
        The complement of $X$ is open and convex. Hence, by Proposition~4 of~\cite{bertsekas_tseng_2006}, we have that $\mathcal{R}(X)=X_\infty$. However, the set $X$ does not satisfy assumption (A1) since, for the direction $d = (0,1) \in X_\infty$
        and any $x$ that lies on the boundary of $X$ (i.e., $x_1^2=x_2$), it is not the case that $x + \alpha d \in X$ for any $\alpha > 0$. Thus, Proposition~12 of~\cite{bertsekas_tseng_2006} cannot be used to claim the existence of solutions in this case. However, if $ \mathcal{K}(f)={\rm Lin}(f)$ (such as, for example, when $f(x_1,x_2)=|x_1|)$ then
        by Theorem~\ref{thm-conf-nonconvx}, the problem $\inf_{x\in X}f(x)$ has a solution.
        \qed
        \end{example}

    \subsection{Constraint Set given by Functional Inequalities}\label{sec-algsets}

    In this section, we consider a problem of the following form: 
    \begin{eqnarray}\label{eq-constrained-func-ineq-prob}
    \hbox{minimize \ \,} && f(x) \cr
    \hbox{subject to } &&  g_j(x) \leq 0, ~j \in \{1,\ldots,m\}, \quad x\in C,
    \end{eqnarray}
     where $C$ is a closed set, and \propf and each $g_j: \ren \to \re \cup \{+\infty\}$ is a proper closed function. 
    Existence of solutions to this problem has been studied in both general  settings \cite{ozdaglar_tseng_2006} as well as in special settings
    \cite{Belousov_Klatte_2002, Eaves-quad-programming, Frank_Wolfe_Theorem_56,   Luo_Zhang_1999, Terlaky-LP}. 
    To the best of our knowledge all the results provided in this section are new, and more general than those in
    \cite{Belousov_Klatte_2002, Eaves-quad-programming, Frank_Wolfe_Theorem_56,   Luo_Zhang_1999, Terlaky-LP,ozdaglar_tseng_2006}, which only provide sufficient conditions.

    We provide necessary and sufficient
    conditions for the existence of solutions to problem \eqref{eq-constrained-func-ineq-prob} based on Theorems~\ref{main-result-2g}--\ref{thm-necessity}. 

    \begin{theorem}[Nonconvex function and nonconvex constraint set]\label{thm-nonconvex-functional}
    Let $C$ be a closed set, and let $f$ and each $g_j$ be proper closed functions with $C\cap\left(\cap_{j=1}^m\dom g_j\right)\cap \dom f \neq \emptyset$.
    Let $X=\{x\in C\mid g_j(x)\le 0, j=1,\ldots,m\}$,
    and assume that $f$ exhibits asymptotically bounded decay on the set $X$ with respect to a coercive function $g$ satisfying the assumptions of Theorem~\ref{main-result-2g}, or  with respect to $g(x)=\|x\|^p$ for some $p\ge0$.
    Then, either one of the following two conditions\\
    \noindent(C1) \quad $C_\infty\cap(\cap_{j=1}^m (L_0(g_j))_\infty)\cap\mathcal{K}(f)\subseteq\mathcal{R}(X)\cap\mathcal{R}(f)$\\
    \noindent (C2) \quad 
    $C_\infty\cap(\cap_{j=1}^m \mathcal{K}(g_j))\cap\mathcal{K}(f)\subseteq\mathcal{R}(X)\cap\mathcal{R}(f)$\\
    is necessary and sufficient for the problem~\eqref{eq-constrained-func-ineq-prob}
    to have a finite optimal value and an optimal solution. 
    \end{theorem}

    \begin{proof}
    By Proposition~\ref{prop-set-intersect}, we have that
    \begin{equation}\label{eq-almost}
    X_\infty\subseteq C_\infty\cap(\cap_{j=1}^m (L_0(g_j))_\infty).
    \end{equation}
    Thus, if condition (C1) holds, then 
    $X_\infty\cap \mathcal{K}(f)\subseteq\mathcal{R}(X)\cap\mathcal{R}(f)$ and
    problem~\eqref{eq-constrained-func-ineq-prob} has a finite optimal value and an optimal solution by Theorem~\ref{main-result-2g} (or Theorem~\ref{main-result}).
    
    If condition (C2) holds, then by Proposition~\ref{prop-level-set} we have
    $(L_0(g_j))_\infty\subseteq\mathcal{K}(g_j)$ for all $j=1,\ldots,m.$
    By combining these relations with~\eqref{eq-almost},
    again we have that $X_\infty\cap \mathcal{K}(f)\subseteq\mathcal{R}(X)\cap\mathcal{R}(f)$ and the result follows as in the preceding case.

    The necessity of these conditions follows by Theorem~\ref{thm-necessity}.
    \end{proof}
    
    In Theorem~\ref{thm-nonconvex-functional}, we could not write the cone $\mathcal{R}(X)$ in terms of such cones of the individual sets defining the set $X$, as there is no particular rule that can be applied here, in general.
    When $X$ and the functions $g_j$ are convex, the conditions (C1) and (C2) of Theorem~\ref{thm-nonconvex-functional} coincide, since $(L_0(g_j))_\infty=\mathcal{K}(g_j)$ for all $j$ by Proposition~\ref{prop-level-set}. 

    We next provide another result for the case when the constraint set $C$ is convex and the functions $g_j$ are convex.
    
    \begin{theorem}[Nonconvex function with convex constraint set]\label{thm-nonconvex-functional1}
    Let assumptions of Theorem~\ref{thm-nonconvex-functional} hold.
    Additionally, assume that the set $C$ is convex and that each $g_j$ is a convex function. Then, either one of the following two conditions\\
    \noindent(C3) \quad 
    $C_\infty\cap(\cap_{j=1}^m \mathcal{K}(g_j))\cap\mathcal{K}(f)\subseteq \mathcal{R}(C)\cap(\cap_{j=1}^m\mathcal{R}(L_0(g_j))\cap\mathcal{R}(f)$\\
    \noindent (C4) \quad 
    $C_\infty\cap(\cap_{j=1}^m \mathcal{K}(g_j))\cap \mathcal{K}(f)\subseteq {\rm Lin}(C)\cap(\cap_{j=1}^m\mathcal{C}(g_j))\cap\mathcal{R}(f)$\\ 
    is necessary and sufficient for problem~\eqref{eq-constrained-func-ineq-prob} to have a finite optimal value and an optimal solution.
    \end{theorem}
    
    \begin{proof} Let the condition (C3) hold.
    Since  $X$ is convex and each $g_j$ is convex, we have that $X_\infty=C_\infty\cap(\cap_{j=1}^m\mathcal{K}(g_j)).$
    Moreover, by
    Proposition~\ref{prop-retractive-intersect} we have 
    \[\mathcal{R}(C)\cap(\cap_{j=1}^m\mathcal{R}(L_0(g_j))\subseteq\mathcal{R}(X).\]
    Thus, the condition (C2) of Theorem~\ref{thm-nonconvex-functional} is satisfied and the result follows.
    Suppose that the condition (C4) holds. Then, since 
    $${\rm Lin}(C)\cap(\cap_{j=1}^m\mathcal{C}(g_j))=
    {\rm Lin}(X)\subseteq\mathcal{R}(X),$$
    the condition (C2) of 
    Theorem~\ref{thm-nonconvex-functional} is satisfied. 
    \end{proof}

Now consider problem~\eqref{eq-constrained-func-ineq-prob} 
where each $g_j$ is a convex polynomial. As a corollary of Theorem~\ref{thm-nonconvex-functional1}, we have the following  result. 
\begin{corollary}[Nonconvex polynomial function with convex polynomial constraints]\label{cor-polyn}
    Consider the problem~\eqref{eq-constrained-func-ineq-prob}, where $C=\re^n$, the objective function $f$ is proper and closed,  and each $g_j$ is a convex polynomial. Then, the problem has a finite optimal value and a solution exists if and only if 
    \[\left(\cap_{j=1}^m{\cal K}(g_j)\right)
    \cap {\cal K}(f)\subseteq \cap_{j=1}^m\mathcal{R}(g_j)\cap \mathcal{R}(f).\]
\end{corollary}

\begin{proof}
    By Example~\ref{ex-lip}, a polynomial of order $p$ asymptotically decays with respect to the function $g(x) = \|x\|^{p}$. 
    Since each $g_j$ is a polynomial, by Lemma~\ref{lem-convpol} we have that 
    $\mathcal{R}(g_j)=\mathcal{C}(g_j)$ for all $j$.
    Therefore, the condition (C4) of Theorem~\ref{thm-nonconvex-functional1} is satisfied, with $C=\re^n$, and the result follows.
\end{proof}

Corollary~\ref{cor-polyn} extends
    Theorem~3 of~\cite{Belousov_Klatte_2002} in two ways. Firstly Theorem~3 applies to the convex case as $f$ is a convex polynomial. Second Theorem~3 provides only sufficient conditions.

The following example 
shows that when the condition 
of Corollary~\ref{cor-polyn} is violated, problem~\eqref{eq-constrained-func-ineq-prob} may not have a solution. 

\begin{example}[Example 2 in~\cite{Luo_Zhang_1999}]\label{luo-zhang-ex2}
    Consider the following problem 
    \begin{equation}\label{LZ-ex-prob}
    \begin{aligned}
    &\hbox{minimize}  &f(x) &= -2x_1x_2 + x_3x_4 + x_1^2 \\
    &\text{subject to \ \ }  &g_1(x) &= x_1^2 - x_3 \leq 0 \\
    &\qquad \ &g_2(x) &= x_2^2 - x_4 \leq 0.
    \end{aligned}
    \end{equation}
The objective polynomial $f$ is nonconvex while both $g_1$ and $g_2$ are convex polynomials. 

By the convexity of the constraint sets, we have 
\[\mathcal{K}(g_1)
=\{(0,d_2,d_3,d_4)\mid d_2\in\re, d_3\ge0, d_4\in\re\},\]
\[\mathcal{K}(g_2)=\{ (d_1,0,d_3,d_4)\mid d_1\in\re, d_3\in\re, d_4\ge0\}.\]
Therefore,
\[\mathcal{K}(g_1)\cap \mathcal{K}(g_2)=\{ (0,0,d_3,d_4)\mid d_3\ge0, d_4\ge0\}.\]
Moreover, by Lemma~\ref{lem-convpol} we have $\mathcal{C}(g_j)=\{0\}$ for $j=1,2$
so that 
$\mathcal{C}(g_1)\cap \mathcal{C}(g_2)=\{0\}.$
The function $f$ is a polynomial of order 2,
so by Lemma~\ref{lem-polynomial} we have 
$\mathcal{K}(f)=\{d\mid f(d)\le 0\}.$
Then, it follows that
\[\mathcal{K}(g_1)\cap \mathcal{K}(g_2)\cap
\mathcal{K}(f)=\{0,0,d_3,d_4)\mid d_3d_4=0, d_3\ge 0, d_4\ge0\}.\]
Since $\mathcal{C}(g_1)\cap \mathcal{C}(g_2)=\{0\}$, we must have
$\mathcal{C}(g_1)\cap \mathcal{C}(g_2)\cap\mathcal{R}(f)=\{0\}.$
Therefore, the condition
\[\mathcal{K}(g_1)\cap \mathcal{K}(g_2)\cap
\mathcal{K}(f)\subseteq \mathcal{C}(g_1)\cap \mathcal{C}(g_2)\cap\mathcal{R}(f)\]
does not hold, since a nonzero direction $d=(0,0,d_3,0)$ with $d_3>0$ belongs to
$\mathcal{K}(g_1)\cap \mathcal{K}(g_2)\cap
\mathcal{K}(f)$ but not to $\mathcal{C}(g_1)\cap \mathcal{C}(g_2)\cap\mathcal{R}(f)$.
However, the optimal value of the problem is $f^*=-1$ which is not attained,
as shown in~\cite{Luo_Zhang_1999}.
\qed
\end{example}

\bmhead{Acknowledgements}
This work has been supported by the ONR award N00014-21-1-2242 and the NSF award CIF 2134256.

\bibliography{manuscript} 

\end{document}